\documentclass[a4paper,12pt]{amsart}
\usepackage{amssymb,amscd}
\usepackage{txfonts}
\usepackage{color}

\newtheorem{thm}{Theorem}[section]
\newtheorem{prop}[thm]{Proposition}
\newtheorem{lem}[thm]{Lemma}
\newtheorem{rem}[thm]{Remark}

\newtheorem{df}[thm]{Definition} 
\newtheorem*{THM*}{Theorem 4.1} 
\numberwithin{equation}{section}
\def\Ker{\operatorname{Ker}}
\def\Im{\operatorname{Im}}
\newcommand{\Q}{\mathbb{Q}}
\newcommand{\Z}{\mathbb{Z}}
\newcommand{\C}{\mathbb{C}}
\newcommand{\F}{\mathbb{F}}
\newcommand{\E}{\mathrm{E}}
\newcommand{\su}{S\!U}

\begin{document}

\title{The Chow rings of the algebraic groups  $\E_{6}$, $\E_{7}$,   and $\E_{8}$}
\author{Shizuo Kaji}
\thanks{The first author is partially supported by Grant-in-Aid for Young Scientists (Start-up)
 20840041, Japan Society for  the Promotion of Science}

\author{Masaki Nakagawa}
\thanks{The second author is partially supported by the Grant-in-Aid for Scientific  Research 
          (C) 21540104, Japan Society for  the Promotion of Science}
\pagestyle{plain}




\address{Department of Applied Mathematics \endgraf 
                        Fukuoka University \endgraf
                        Fukuoka 814-0180   \\ Japan}
\email{kaji@math.sci.fukuoka-u.ac.jp} 

\address{Department of General Education \endgraf
                      Kagawa National College of Technology \endgraf
                      Takamatsu 761-8058 \\ Japan}
\email{nakagawa@t.kagawa-nct.ac.jp} 

\begin{abstract}
     We determine the Chow rings of the complex algebraic groups of the exceptional type $\E_6, \E_7$,  and $\E_8$,
     giving  the explicit generators represented by the pull-back images of Schubert varieties 
     of the corresponding flag varieties. This is a continuation of the work of R. Marlin on
     the computation of the Chow rings of $\mathrm{SO}_n$, $\mathrm{Spin}_n$, $\mathrm{G}_2$, 
     and $\mathrm{F}_4$.   
     Our method is based on Schubert calculus of the corresponding flag varieties, which has its own interest.
\end{abstract} 

\maketitle

\section{Introduction}
The problem of computing the Chow rings of complex algebraic groups
dates back to the paper by Grothendieck \cite{Gro58}, where he showed the Chow rings of 
$\mathrm{SL}_{n}$ and $\mathrm{Sp}_{2n}$ are trivial.
Later R. Marlin computed the Chow rings of $\mathrm{SO}_{n}$, $\mathrm{Spin}_{n}$,
 $\mathrm{G}_{2}$, and $\mathrm{F}_{4}$ (\cite{Mar74}, \cite{Mar74(2)}).
The purpose of the present article is to give the explicit descriptions of the Chow rings 
 for the remaining cases of $\E_{6}$, $\E_7$,  and  $\E_{8}$\footnote{Recently,
 H. Duan and X. Zhao also determined the Chow rings of $\E_{6}$, $\E_{7}$, and $\E_{8}$ by 
 a different approach from ours (\cite{DZ07}).}.

Let $G$ be a  simply-connected simple  algebraic group over the field of complex numbers $\C$,
$B$ a fixed  Borel subgroup  of $G$, and $H$ a maximal (algebraic) torus contained in $B$. 
The homogeneous space $G/B$ is known to be a projective variety,  called the  {\it flag variety}.
The well-known basis theorem \cite{Che58} tells us  that the Chow ring $A(G/B)$ of $G/B$ has a
distinguished $\Z$-module basis which consists of {\it Schubert varieties}.
Our approach is to describe $A(G)$ via pull-back images of Schubert varieties
 through the quotient map $p: G \to G/B$.  To be precise, the mothod proceeds as follows: 
 Let 
\begin{equation}  \label{eqn:ch.hom}
  c_{G}:  S(\hat{H}) \longrightarrow  A(G/B)  
\end{equation} 
be the characteristic homomorphism for $G$, where  $S(\hat{H})$ denotes 
the symmetric algebra over $\Z$ of  the character group $\hat{H}$ of $H$  (\cite[(4.1)]{Gro58}).   
Then Grothendieck showed (\cite[p.21, R{\scriptsize  EMARQUES} $2^{\circ}$]{Gro58}, 
see also \cite[\S 3]{Bri98} for the proof)   
that $p^{*}: A(G/B) \longrightarrow A(G)$ is surjective and its kernel is the ideal generated 
by $c_{G}(\hat{H})$, in other words, the divisor classes on $G/B$. Therefore the determination 
of $A(G)$ reduces to that of $A(G/B)$. 
In principle, the ring structure of $A(G/B)$ could  be determined 
by the classical Chevalley formula \cite{Che58}, but in practice,  it could  not be applicable 
to the cases of higher rank groups.
On the other hand, the Chow ring  $A(G/B)$   is  isomorphic  to  the integral cohomology ring 
$H^*(G/B;\Z)$    via the cycle map  which assigns to each Schubert variety its fundamental class 
(\cite[\S 6]{Gro58}, \cite[Example 19.1.11]{Ful98}).
Hence we work with the cohomology instead of the Chow ring,
which allows us to access the computational facility of algebraic topology.

To describe the cohomology ring $H^*(G/B;\Z)$, we have two different ways, namely,
the {\it Schubert presentation} and the {\it Borel presentation}.
In the Schubert presentation, an additive basis for $H^*(G/B;\Z)$ is given by 
Schubert classes corresponding to the Schubert varieties (\cite{BGG73},  \cite{Bor54}, \cite{Che58}). 
However, the multiplicative structure  among Schubert classes  is difficult to determine, and hence  
 it does not seem to fit our purpose of computing $A(G)$.   
 In fact, Marlin's computation of $A(G)$ for $G = \mathrm{SO}_{n}$,  $\mathrm{Spin}_{n}$, $\mathrm{G}_{2}$, and 
$\mathrm{F_4}$  in \cite{Mar74}, where he considered the Schubert presentation only, 
becomes unmanageable when applied to 
the remaining exceptional groups $\E_{6}, \E_{7}$, and $\E_{8}$.
In \cite{Bor53}, A. Borel gives another description for $H^*(G/B;\Z)$ in terms of
the ring of invariants of the Weyl group action on   $S(\hat{H})$, which we call the Borel presentation today. 
This presentation has an advantage that  the ring structure of $H^*(G/B;\Z)$ is relatively
easy to see.  However,  the generators in this presentation have little geometric meaning.  
There is a connection between these two presentations discovered by Bernstein-Gelfand-Gelfand  \cite{BGG73} 
and Demazure \cite{Dem73}.
More specifically,  they introduced  a series of   operators called 
the  {\it divided difference operators}  acting on  $S(\hat{H})$.
Using these operators, we can  express the ring generators of 
$H^*(G/B;\Z)$ obtained by the Borel presentation in terms of Schubert classes.
Computational difficulties in the cases of the  higher rank exceptional groups $G= \E_{l} \ (l = 6,7,8)$
 are overcome by  an appropriate choice of the generators for $H^*(\E_{l}/B;\Z)$
discovered by Toda \cite{Toda75}  and Toda-Watanabe \cite{Toda-Wat74} (see \S \ref{subsec:toda-watanabe-basis}).
Once we obtain the correspondence between the Borel and Schubert presentations,
 the Chow rings $A(\E_{l})$  are determined immediately using  the  
result of Grothendieck mentioned earlier.   In fact, the second author  simplified
Marlin's computation by  making use of the same method as in  this paper  (\cite{Nak08}).

 Then our main result is   stated as follows: 
\begin{thm}      \label{thm:main_theorem} 
For $G = \E_{l} \ (l=6,7,8)$,  we denote by $p:G \longrightarrow G/B$ the natural  projection, by 
$p^* : A(G/B) \longrightarrow A(G)$ the induced pull-back homomorphism,  by $w_0$ 
the longest element of the Weyl group of $G$, and  by $s_{i}$ the reflection corresponding to 
the simple root  $\alpha_{i} \ (1 \leq i \leq l)$ $($for the notation, see $\S \ref{subsec:Schubert}$ and 
$\S \ref{subsec:toda-watanabe-basis})$.

{\rm (1)} The Chow ring of $\E_{6}$ is given by
  \[  A(\E_{6}) = \Z[X_{3}, X_{4}]/(2X_{3}, 3X_{4}, X_{3}^2, X_{4}^3),  \]
  where $X_{3}$ and $X_{4}$ are the images under $p^*$ of the elements of $A(\E_{6}/B)$  defined by 
  the Schubert varieties  $X_{w_0s_{5}s_{4}s_{2}}$ and $X_{w_0s_{6}s_{5}s_{4}s_{2}}$ respectively.
   
{\rm (2)} The Chow ring of $\E_{7}$ is given by
      \[    A(\E_{7})  =  \  \Z[X_{3}, X_{4}, X_{5}, X_{9}]  
                            /(2X_{3}, 3X_{4}, 2X_{5}, X_{3}^2, 2X_{9}, X_{5}^2, X_{4}^3, X_{9}^2),   \]   
  where $X_{3}, X_{4}, X_{5}$,  and $X_{9}$ are the images under $p^*$ of the elements of  $A(\E_{7}/B)$ 
  defined by the Schubert varieties  
    $X_{w_0s_{5}s_{4}s_{2}}, 
     X_{w_0s_{6}s_{5}s_{4}s_{2}}, 
     X_{w_0s_{7}s_{6}s_{5}s_{4}s_{2}}$,  
and $X_{w_0s_{6}s_{5}s_{4}s_{3}s_{7}s_{6}s_{5}s_{4}s_{2}}$
  respectively. 
  
 {\rm (3)} The Chow ring of $\E_{8}$ is given by
  \begin{align*}
        A(\E_{8})  = & \Z[X_{3}, X_{4}, X_{5}, X_{6}, X_{9}, X_{10}, X_{15}]  \\
   &      /(2X_{3}, 3X_{4}, 2X_{5}, 5X_6,  2X_{9}, 3X_{10}, X_{4}^3, 2X_{15},  X_{9}^2, X_{5}^4, X_{3}^8, X_6^5, X_{10}^3, X_{15}^2), 
      \end{align*}
  where $X_{3}, X_{4}, X_{5}, X_6, X_{9}, X_{10}$, and $X_{15}$ are the images under $p^*$ of the elements of  $A(\E_{8}/B)$ 
  defined by the Schubert varieties  
   $X_{w_0s_{5}s_{4}s_{2}}, 
    X_{w_0s_{6}s_{5}s_{4}s_{2}}, 
    X_{w_0s_{7}s_{6}s_{5}s_{4}s_{2}}, 
    X_{w_0s_{1}s_{3}s_{6}s_{5}s_{4}s_{2}},
    X_{w_0s_{1}s_{5}s_{4}s_{3}s_{7}s_{6}s_{5}s_{4}s_{2}}$,
   $X_{w_0s_{1}s_{6}s_{5}s_{4}s_{3}s_{7}s_{6}s_{5}s_{4}s_{2}}$, and
   $X_{w_0s_{1}s_{3}s_{4}s_{2}s_{7}s_{6}s_{5}s_{4}s_{3}s_{8}s_{7}s_{6}s_{5}s_{4}s_{2}}$  respectively.  
\end{thm}

Combining Theorem \ref{thm:main_theorem} with 
 the results from \cite{Mar74}, \cite{Mar74(2)}, we have the following table:

{ \renewcommand\arraystretch{2.8}
\begin{tabular}{lcll}
$G$ & $A(G)$ & generators \\
\hline
$\mathrm{SL}_{n} $   & $\Z$ \\
$\mathrm{Sp}_{2n}$   & $\Z$  \\
$\mathrm{SO}_{2n}$   
& $\dfrac{\Z[X_{1}, X_{3}, X_{5}, \ldots, X_{2[\frac{n}{2}]-1}]}{(2X_{i}, \; X_{i}^{p_{i}})},  \;  
  p_{i} = 2^{[\log_{2}\frac{n-1}{i}] + 1}$ 
& $[n],[n-2,n], \ldots, [1, \ldots, n-3, n-2, n]$ \\
$\mathrm{Spin}_{2n}$ 
& $\dfrac{\Z[X_{3}, X_{5}, \ldots, X_{2[\frac{n}{2}]-1}]}{(2X_{i}, \; X_{i}^{p_{i}})},  \;  
p_{i} = 2^{[\log_{2}\frac{n-1}{i}] + 1}$ 
& $[n-2, n], \ldots, [1, \ldots, n-3, n-2, n]$ \\
$\mathrm{SO}_{2n+1}$ 
& $\dfrac{\Z[X_{1}, X_{3}, X_{5}, \ldots, X_{2[\frac{n+1}{2}]-1}]}{(2X_{i},\; X_{i}^{p_{i}})}, \;  
p_{i} = 2^{[\log_{2}\frac{n}{i}] + 1}$ 
& $[n], [n-1, n], \ldots, [1, \ldots, n]$   \\
$\mathrm{Spin}_{2n+1}$ 
& $\dfrac{\Z[ X_{3}, X_{5}, \ldots, X_{2[\frac{n+1}{2}]-1}]}{(2X_{i}, \; X_{i}^{p_{i}})}, \;  
p_{i} = 2^{[\log_{2}\frac{n}{i}] + 1}$ 
& $[n-1, n], \ldots, [1, \ldots, n]$ \\
$\mathrm{G}_{2}$ & $\dfrac{\Z[X_{3}]}{(2X_{3}, X_{3}^2)}$ & $[121]$ \\
$\mathrm{F}_{4}$ & $\dfrac{\Z[X_{3}, X_{4}]}{(2X_{3}, 3X_{4}, X_{3}^2, X_{4}^3)}$ & $[123], [1234]$ \\
$\E_{6}$ & $\dfrac{\Z[X_{3}, X_{4}]}{(2X_{3}, 3X_{4}, X_{3}^2, X_{4}^3)}$ & $[542], [6542]$ \\
$\E_{7}$ & $\dfrac{\Z[X_{3}, X_{4}, X_{5}, X_{9}]}{(2X_{3}, 3X_{4}, 2X_{5}, X_{3}^2, 2X_{9}, X_{5}^2, X_{4}^3, X_{9}^2)}$ 
& \parbox{5cm}{$[542], [6542]$,\\ $[76542], [654376542]$}\\
$\E_{8}$ & $\dfrac{\Z[X_{3}, X_{4}, X_{5}, X_{6}, X_{9}, X_{10}, X_{15}]}
                {  \parbox{6cm}{$ (2X_{3}, 3X_{4}, 2X_{5}, 5X_6,  2X_{9}, 3X_{10}, X_{4}^3, $\\
                   $\qquad 2X_{15},  X_{9}^2, X_{5}^4, X_{3}^8, X_6^5, X_{10}^3, X_{15}^2)$}}$ & 
                   \parbox{5cm}{$[542], [6542]$,\\ $[76542], [136542], [154376542]$, \\
                   $[1654376542], [134276543876542]$
}\\
\end{tabular}}
Here  a Schubert variety $X_{w_{0}s_{i_1}s_{i_2}\ldots}$ is abbreviated to  $[i_{1} i_{2} \ldots]$ 
in the third column of the table.
\begin{rem}\rm
We note that the mod $p$ Chow rings $A(G;\F_{p}) := A(G) \otimes_{\Z}  \F_{p}$ 
were studied by V. Kac in \cite{Kac85}, where he showed that $A(G;\F_{p})$ is isomorphic to $\F_{p}$ 
for a non-torsion prime $p$ of $G$,  whereas for a torsion prime $p$, $A(G;\F_{p})$ is 
the ``polynomial part'' of $H^*(G;\F_{p})$. 
\end{rem}

The organization of this paper is as follows: 
In \S \ref{sec:cohomology}, we review basic facts on Schubert calculus and fix our notations.
In \S \ref{sec:integral_cohomology},  we recall the Borel presentation of the integral cohomology rings of 
$H^*(G/B;\Z)$ for $G = \E_{l} \;  (l = 6, 7, 8)$ and  convert those results to the Schubert presentation using the 
divided difference operators. Finally in \S \ref{sec:Chow}, we give descriptions of 
the Chow rings $A(G)$ for  $G = \E_{l} \;  (l = 6, 7, 8) $.

\textbf{Acknowledgments.}  
We thank  Nobuaki Yagita for giving us useful comments concerning the Chow rings of algebraic groups.   
He  pointed out that  the ring  structure of $A(G)$  can  also be obtained  by using the technique of \cite{Yag05}.   
We also thank  Haibao  Duan and Xuezi  Zhao for discussion about the cohomology of flag varieties. 
Finally, we thank Mamoru Mimura for giving us various suggestions.


\section{Cohomology of flag varieties} \label{sec:cohomology}
In the area of Schubert calculus, we frequently use
two ways of describing the integral cohomology of flag varieties;
the Schubert presentation and the Borel presentation.
The former gives a geometric basis for the cohomology whose product structure is hard to know,
while the latter allows   us purely algebraic treatment of the cohomology. 
For our purpose, we need  both presentations.

\subsection{}  \label{subsec:Schubert} 
We begin with a brief review of the Schubert presentation  of the integral cohomology rings 
of  flag varieties (\cite{BGG73}, \cite{Bor54},  \cite{Che58}).  

Let us denote by  $\Delta$  the root system  of $G$ relative to  $H$,  by  $\Pi$  the system of  simple roots,
by $N_{G}(H)$ the normalizer of $H$ in $G$. Then the Weyl group $W$ of $G$ is defined by $N_{G}(H)/H$.
Denote  by  $s_{\alpha}$ the simple reflection, i.e., the  reflection corresponding to the simple root 
$\alpha \in \Pi$, and by  $S = \{ s_{\alpha} \, | \, \alpha \in \Pi \}$ the set of simple reflections. 
Then it is known that the  Weyl group $W$ of $G$ is generated by the simple reflections. 
Denote by $l(w)$ the length of  an element $w \in W$ with  respect to $S$.

As is well known (see for example \cite[14.12]{Bor91}), $G$ has the Bruhat decomposition 
 \begin{equation*} 
      G  =  \coprod_{w \in W} B \dot{w} B,   
  \end{equation*} 
where $\dot{w}$ denotes any   representative of  $w \in W$.  It induces a cell  decomposition
\begin{equation*} 
      G/B =  \coprod_{w \in W}  B \dot{w} B/B,   
 \end{equation*} 
where $X_{w}^{\circ} = B \dot{w} B/B \cong \C^{l(w)}$ is called the {\it Schubert cell}.
The {\it Schubert variety}  $X_{w}$ is defined  to be the closure $\overline{X_{w}^{\circ}}$ of $X_{w}^{\circ}$. 
The fundamental class $[ X_w ]$  of $X_{w}$ lies in $H_{2l(w)}(G/B;\Z)$.  We define a cohomology class 
$Z_w \in H^{2l(w)}(G/B;\Z)$ as the Poincar\'{e} dual of $[X_{w_0w}]$, where $w_0$ is the longest 
element of $W$.  We call $Z_w$ the {\it Schubert class}  corresponding to $w \in W$.  The Schubert classes 
$\{ Z_{w} \}_{w \in W}$ form an additive basis for  the free $\Z$-module $H^*(G/B;\Z)$;  we refer to $\{ Z_{w} \}_{w \in W}$ 
as the {\it Schubert basis}. 

The product of two Schubert classes can be expressed as a $\Z$-linear combination of the Schubert basis. 
In order to complete  the multiplicative structure of $H^*(G/B;\Z)$, we have to compute 
the {\it structure constants} $a_{u, v}^{w}$ for $u, v, w \in W$. These integers $a_{u, v}^{w}$ 
are  defined by the following equation:
\begin{equation*} 
       Z_{u} \cdot Z_{v} = \hspace{-0.8cm} 
                             \sum_{\tiny{ \begin{array}{ccc}
                                           & w \in W, \\ 
                                           &  l(u) + l(v) = l(w)
                                        \end{array} 
                                       }
                                  } \hspace{-0.8cm}  a_{u,v}^{w} \, Z_{w}.    
\end{equation*} 
Computing the structure constants $a_{u, v}^{w}$ becomes  a hard task in general.   
Indeed, one of the main problems  of {\it Schubert calculus} is to give 
a formula for structure constants (see for example \cite{Pra05}).  

\begin{rem} \rm
 In \cite{Duan05}, H. Duan find an effective algorithm for computing the structure constants.
Based on it, he and X. Zhao gave  descriptions of  the integral cohomology rings
of flag varieties in terms of Schubert classes $($\cite{DZ08}$)$.  
\end{rem}

To a subset $I_{P}$ of simple roots $\Pi$, we associate the parabolic subgroup $P \; (\supset B)$ of $G$
whose Weyl group\footnote{More precisely, $W_{P}$ is the Weyl group of the reductive part of $P$.} $W_{P}$ is 
the subgroup of $W$ generated by $\{ s_{\alpha}  \mid \alpha  \in I_{P} \}$.
We can take the explicit minimal coset representatives $W^{P}$ of $W/W_{P}$ 
 (see \cite[\S 1.10]{Hum90}) as 
   \[ W^P  = \{ w \in W \mid   l(w s_{\alpha}) =   l(w) + 1 \;  \text{for all}   \; \alpha \in \;  I_P \}. \]
For $w \in W^P$, the Schubert class $Z_{w} \in H^{2l(w)}(G/P;\Z)$ of the partial flag variety $G/P$ is defined as the class
dual to   $[\overline{B w_{0} w  P/P}]$. 
Similarly to the case of $G/B$, the Schubert classes 
$\{ Z_{w} \}_{w \in W^{P}}$ form an additive basis for  the free $\Z$-module $H^*(G/P;\Z)$.  
Consider  the following  fibration: 
\begin{equation}\label{eqn:bundle}
   P/B  \overset{i}{\hookrightarrow}    G/B \xrightarrow{\pi} G/P.
 \end{equation}
Since the fibre $P/B$ and the base $G/P$ have even dimensional cohomology,  
the  Serre spectral sequence  with integer coefficients for the fibration (\ref{eqn:bundle}) 
collapses at the $E_{2}$-term. 
Therefore the projection $\pi$ induces an  inclusion $\pi^{*}:  H^*(G/P;\Z) \hookrightarrow  H^*(G/B;\Z)$
which is compatible with  the inclusion $W^P \subset W$, the inclusion $i$ induces a surjection 
$i^{*}: H^*(G/B;\Z) \twoheadrightarrow H^*(P/B;\Z)$, and $\Ker i^*$ is an ideal of $H^*(G/B;\Z)$ generated by 
$\pi^{*}H^{+}(G/P;\Z)$, where $H^{+}(G/P;\Z) = \bigoplus_{i > 0}H^{i}(G/P;\Z)$.   
Furthermore, it is known that $\Im \pi^{*} = H^{*}(G/P;\Z)  \subset H^{*}(G/B;\Z)$ coincides with 
the set of $W_{P}$-invariant elements of $H^{*}(G/B;\Z)$ (\cite[Theorem 5.5]{BGG73}). 
From this, we can describe the ring structure of $H^*(G/B;\Z)$ in terms of $H^*(G/P;\Z)$ and $H^*(P/B;\Z)$
(cf. \cite[Lemma 1,1]{Toda-Wat74},  \cite[\S 3]{Mar74}, \cite[\S 2.8]{DZ08}). 
Certain choice of $P$ reduces the computation of $H^*(G/B;\Z)$ greatly via the ``splitting'' above. 
Especially for our purpose of determining the Chow ring $A(G) \cong H^*(G/B;\Z)/ \left( H^2(G/B;\Z) \right)$,
we will choose $P$ so that $H^*(P/B;\Z)$ is generated by degree two elements (see \S \ref{subsec:toda-watanabe-basis}).

\subsection{}    \label{subsec:Borel}   
  In this subsection, we review the  Borel  presentation of  the cohomology  of flag  varieties  (\cite{Bor53}).

Let $K$ be a maximal compact subgroup of $G$ and $T = K \cap H$ a compact maximal  torus of $K$. 
Then the inclusion $K \hookrightarrow G$ induces a diffeomorphism  $K/T  \cong G/B$.
The inclusion   $T \hookrightarrow K$ induces a  fibration 
  \[   K/T \overset{\iota}{\longrightarrow}  BT \overset{\rho}{\longrightarrow} BK,    \]
where $BT$ (resp. $BK$) denotes the classifying space of $T$ (resp. $K$). 
The induced homomorphism in cohomology   
   \begin{equation} \label{eqn:ch.hom2}
      c  = \iota^*:  H^*(BT;\Z) \longrightarrow H^*(K/T;\Z)   
   \end{equation} 
is called  {\it Borel's characteristic homomorphism}\footnote{It is not difficult to see that the symmetric algebra 
$S(\hat{H})$ is isomorphic to $H^*(BT;\Z)$, and Borel's characteristic homomorphism is identified with 
the  characteristic homomorphism   $c_{G}:S(\hat{H}) \longrightarrow A(G/B)$.}.  The Weyl group $W$ of $K$ acts naturally on $T$, 
hence on $H^{2}(BT;\Z)$.   This  action of $W$  extends to the whole  $H^*(BT;\Z)$ and also to 
$H^*(BT;\F)= H^*(BT;\Z) \otimes_{\Z} \F$, where $\F$ is any field. 
Then one of Borel's results can be stated as follows: 

\begin{thm}[Borel \cite{Bor53}]  \label{thm:Borel}
  Let $\F$ be a field of characteristic zero. Then Borel's characteristic homomorphism induces 
  an isomorphism  
     \[  \overline{c}: H^*(BT;\F)/(H^{+}(BT;\F)^{W}) \longrightarrow H^*(K/T;\F),  \]
where $(H^{+}(BT;\F)^{W})$ is the ideal of $H^*(BT;\F)$ generated by the $W$-invariants of 
positive degrees. 
\end{thm} 

In particular,  the rational cohomology ring $H^*(K/T;\Q)$ is easy to describe. 
 For the classical types, the rings of invariants of the Weyl groups with rational 
coefficients are well-known.  For the exceptional types, see for example \cite{Meh88}.      
In order to determine the integral cohomology 
ring $H^*(K/T;\Z)$, we need further considerations.  In \cite{Toda75}, Toda established a method to 
describe the integral cohomology ring $H^*(K/T;\Z)$ by a minimal system of generators and relations, 
from the mod $p$ cohomology rings $H^*(K;\Z/p\Z)$ for all primes $p$ and the rational cohomology ring $H^*(K/T;\Q)$. 
Along the line  of Toda's method,  the integral cohomology rings of flag varieties for 
$G = \E_{l} \ (l = 6,7, 8)$
have been concretely determined  by Toda-Watanabe \cite{Toda-Wat74} and the second author \cite{Nak01}, \cite{Nak09.2}.
We emphasize that the ring structure is completely determined there.  
Throughout the rest of the paper, we fix maximal compact subgroups of  $\E_{l} \; (l = 6, 7, 8)$
and denote them by  $E_{l} \; (l = 6, 7, 8)$ respectively.

 \section{Integral cohomology ring of $E_l/T \; (l = 6, 7, 8)$  }   \label{sec:integral_cohomology}  
 In this section, we focus on the cases of the exceptional groups $E_l \ (l=6,7,8)$.

\subsection{}  \label{subsec:toda-watanabe-basis}
First we introduce a handy set of generators of $H^2(BT;\Z)$ for $K=E_l \ (l=6,7,8)$.
 Following  \cite{Bour68}, we take the simple roots 
  $\Pi  = \{ \alpha_{i} \}_{1 \leq i \leq l}$.  For instance, the Dynkin diagram of $E_{8}$ is given as follows: 
\setlength{\unitlength}{1mm}
\begin{figure}[ht]
      \begin{picture}(60,30)
         \multiput(-15,20)(15,0){7}{\circle{2}} 
         \put(15,5){\circle{2}}     
         \multiput(-14,20)(15,0){6}{\line(1,0){13}}
         \put(15,19){\line(0,-1){13}}          
         \put(-15,27){\makebox(0,0)[t]{$\alpha_{1}$}}    
         \put(0,27){\makebox(0,0)[t]{$\alpha_{3}$}}
         \put(15,27){\makebox(0,0)[t]{$\alpha_{4}$}}
         \put(30,27){\makebox(0,0)[t]{$\alpha_{5}$}}
         \put(45,27){\makebox(0,0)[t]{$\alpha_{6}$}}
         \put(60,27){\makebox(0,0)[t]{$\alpha_{7}$}}
         \put(75,27){\makebox(0,0)[t]{$\alpha_{8}$}}
         \put(23,5){\makebox(0,0)[r]{$\alpha_{2}$}}    
    \end{picture}
      \caption{Dynkin diagram of $E_{8}$}\label{figure:dynkin}
\end{figure}
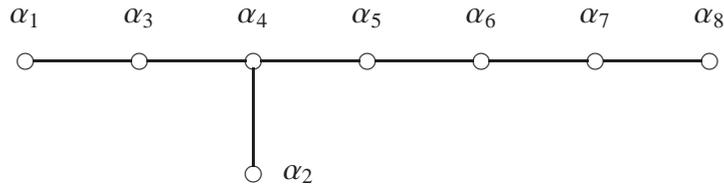

We denote by  $\{ \omega_{i} \}_{1 \leq i \leq l}$ the  corresponding 
fundamental weights. As is customary,  we regard roots and weights  as elements of $H^{2}(BT;\Z)$.  
Let $s_i \ (1 \leq i \leq l)$ denote the  reflection corresponding to the simple root  
$\alpha_{i} \ (1 \leq i \leq l)$. Then the Weyl group  $W(E_{l})$ of $E_{l}$ is 
generated by $s_{i} \ (1 \leq i \leq l)$.  

  Following  \cite[\S 4]{Toda-Wat74}, \cite[\S 1]{Wat75},  and \cite[\S 2]{Nak09},  we set 
    \begin{equation}  \label{eqn:t_i}
        \begin{array}{llll} 
          t_{1} = & \hspace{-0.2cm} -\omega_{1} + \omega_{2}, \;     
          t_{2} = \omega_{1} + \omega_{2} - \omega_{3}, \;    
          t_{3} = \omega_{2} + \omega_{3} - \omega_{4}, \medskip \\
          t_{i} = & \omega_{i} - \omega_{i+1} \; (4 \le i < l),  \; 
          t_{l} = \omega_{l}, \; 
              t = \omega_{2}, \; 
          c_{i} = e_{i}(t_{1}, \ldots, t_{l}) \;    (1 \leq i \leq l),  \medskip   
        \end{array}   
    \end{equation} 
where $e_{i}(t_{1}, \ldots, t_{l})$ denotes the $i$-th elementary symmetric polynomial in the variables 
$t_{1}, \ldots, t_{l}$. In this setting, we have 
    \begin{align*} 
          H^*(BT;\Z) & = \Z[\omega_{1}, \omega_{2}, \ldots, \omega_{l} ]  \\
                     & = \Z[t_{1}, t_{2}, \ldots, t_{l}, t]/(c_{1} - 3t).     
    \end{align*} 
 Since we assume  the simply-connectedness  of the groups, Borel's characteristic homomorphism
 restricted in degree $2$ is an isomorphism $H^{2}(BT;\Z)  \cong    H^{2}(E_{l}/T;\Z)$.
Under this isomorphism, we denote the images of $t_{i} \ (1 \leq i \leq l)$ and  $t$  by the same symbols. 
Thus $H^{2}(E_{l}/T;\Z)$ is a free $\Z$-module generated by $t_{i} \ (1 \leq i \leq l)$ and $t$
with a relation $c_{1} = 3t$. 

As for the action of the Weyl group on $H^{2}(BT;\Z)$ ($\cong H^{2}(E_{l}/T;\Z)$), one can easily see  that 
   \begin{itemize} 
           \item [(i)] $s_{i} \ (i \neq 2)$ acts  on $\{ t_{i} \}_{1 \leq i \leq l}$
                       as the transpositions and trivially on $t$.  
           \item [(ii)] The action of $s_{2}$ on $\{ t_{i} \}_{1 \leq i \leq l}$  
                        and $t$  is given by  
                        \begin{equation*}  
                            \begin{array}{lll} 
                                s_{2}(t_{i}) & = \left \{\hspace{-0.3cm}  
                                        \begin{array}{lll} 
                                             & t - t_{1} - t_{2} - t_{3} + t_{i} \quad \text{for} \quad 
                                               1 \leq i \leq  3,  \medskip \\
                                             & t_{i}  \quad  \text{for} \quad  4 \leq i \leq  l, \medskip  
                                        \end{array}   
                                                 \right. \\
                                s_{2}(t)   &=  2t - t_{1} - t_{2} - t_{3}. \medskip 
                        \end{array}
                       \end{equation*}           
        \end{itemize}

Let $P_{2}$ be the centralizer of a one dimensional torus determined by $\alpha_{i} = 0 \; (i \neq 2)$
\footnote{The subgroup $P_{2}$ corresponds to the maximal parabolic subgroup of $\E_{l}$ associated to the 
subset $\Pi \setminus \{ \alpha_{2} \}$.}. 
Then  the  Weyl group $W_{P_{2}}$  of $P_{2}$  is generated by $s_{i} \ (i \neq 2)$.
Notice that the elements $t_{i} \; (1 \leq i \leq l)$ and $t$  in $H^{2}(E_{l}/T;\Z)$ 
belong to  an   orbit under the action of  $W_{P_{2}}$.  
This technical choice of a ``maximal parabolic subgroup'' of $E_{l}$  is made by paying attention to the symmetry seen 
from the Dynkin diagram of type $E$   exceptional groups (see Figure \ref{figure:dynkin});  
Disregarding the root $\alpha_{2}$, the Dynkin diagram of type $E_{l}$ is the same as that of type $A_{l-1}$.
Therefore, by  (\ref{eqn:bundle}),  we have the following bundle: 
\[
\su(l)/T' \cong P_{2}/T  \overset{i}{\hookrightarrow}   E_{l}/T   \overset{\pi}{\to}   E_{l}/P_{2}, 
\]
where $T'$ denotes the standard maximal torus of $\su(l)$.  As mentioned in \S \ref{subsec:Schubert},  
the Serre spectral sequence  with integer coefficients for the above fibration collapses at the $E_{2}$-term.  
From this and the fact that  $H^*(\su(l)/T';\Z)$ is generated by  degree two classes, 
it is not hard to see that $H^*(E_{l}/T;\Z)$ is  generated  multiplicatively  by $H^{2}(E_{l}/T;\Z)$ and $\Im \pi^{*} = 
H^{*}(E_{l}/P_{2};\Z)$. Especially, ring generators of degrees greater than two can be chosen so that 
they lie in $H^{*}(E_{l}/P_{2};\Z) \subset H^*(E_{l}/T;\Z)$ (see Theorems \ref{thm:E_6/T}, \ref{thm:E_7/T}, and \ref{thm:E_8/T} below).  Hence   any class in 
$H^*(E_l/T;\Z)/ \left( H^2(E_l/T;\Z) \right)$ $( \cong A(\E_{l}) )$ can be represented by an element 
of $H^*(E_l/P_2;\Z) \subset H^*(E_l/T;\Z)$.
This is the crucial point for the computation in the sequel.

 \subsection{}   \label{subsec:E_l/T}
Using the basis for $H^{2}(E_{l}/T;\Z)$ described in the previous subsection,
we give  concrete descriptions of the  integral cohomology rings of $E_{l}/T$,
following the results of \cite{Toda-Wat74}, \cite{Nak01},  and \cite{Nak09.2}.

\begin{thm}[Toda-Watanabe \cite{Toda-Wat74}, Theorem B] \label{thm:E_6/T}
  The integral cohomology ring of $E_{6}/T$ is 
  \[ H^{*}(E_{6}/T;\Z) = \Z [t_{1},  \ldots, t_{6}, t, \gamma_{3}, \gamma_{4} ]
              /(\rho_{1}, \rho_{2}, \rho_{3}, \rho_{4}, \rho_{5}, \rho_{6}, \rho_{8},  \rho_{9}, \rho_{12}),       \]
   where
  \begin{align*}
     \rho_{1}  &= c_{1} - 3t, \;   
     \rho_{2}   = c_{2} - 4t^2, \;    
     \rho_{3}   = c_{3} - 2\gamma_{3}, \;  
     \rho_{4}   = c_{4} + 2t^4 - 3\gamma_{4}, \\ 
     \rho_{5}  &= c_{5} - 3t\gamma_{4} + 2t^2\gamma_{3}, \; 
     \rho_{6}   = \gamma_{3}^2 + 2c_{6} - 3t^2 \gamma_{4} + t^6, \\ 
     \rho_{8}  &= 3 \gamma_{4}^2 - 6t\gamma_{3}\gamma_{4} - 9t^2 c_{6} + 15t^4 \gamma_{4} - 6t^5 \gamma_{3} - t^8, \; 
     \rho_{9}   = 2c_{6}\gamma_{3} - 3t^3 c_{6},  \medskip \\
     \rho_{12} &= 3c_{6}^2 - 2 \gamma_{4}^3 + 6t\gamma_{3} \gamma_{4}^2 + 3t^2c_{6}\gamma_{4} 
                  + 5t^3 c_{6}\gamma_{3}- 15t^4 \gamma_{4}^2 - 10t^6c_{6} 
                  + 19t^8 \gamma_{4} - 6t^9 \gamma_{3} - 2t^{12}. \medskip 
\end{align*}
\end{thm}
\begin{rem} \rm  
This description is slightly different from the one in \cite{Toda-Wat74}.
Toda-Watanabe \cite{Toda-Wat74} make 
use of the integral cohomology ring of the Hermitian symmetric space $E III = E_{6}/T^1 \! \cdot \! \mathit{Spin}(10)$
to determine the higher relations $\rho_{9}$ and $\rho_{12}$.
For our purpose, it is convenient to have  a description directly related to the ring of invariants of the Weyl group 
$W(E_{6})$.  We will give the outline of the computation  in \S \ref{Borel_presentation}  for  convenience of the reader. 
\end{rem}

\begin{thm} [Nakagawa \cite{Nak01}, Theorem 5.9]  \label{thm:E_7/T} 
The integral cohomology ring of $E_{7}/T$ is 
     \begin{align*} 
             H^{*}(E_{7}/T;\Z) =\  &  \Z [t_{1}, \ldots, t_{7}, t, 
                                            \gamma_{3}, \gamma_{4}, \gamma_{5}, \gamma_{9} ]  \\
                                   &  /(\rho_{1}, \rho_{2}, \rho_{3}, \rho_{4}, \rho_{5}, \rho_{6}, \rho_{8}, \rho_{9}, 
                                        \rho_{10}, \rho_{12}, \rho_{14}, \rho_{18}),   
    \end{align*}
  where
  \begin{align*} 
     \rho_{1}  &= c_{1} - 3t, \;    
     \rho_{2}   = c_{2} - 4t^2, \;    
     \rho_{3}   = c_{3} - 2\gamma_{3}, \;   
     \rho_{4}   = c_{4} + 2t^4 - 3\gamma_{4}, \\
     \rho_{5}  &= c_{5} - 3t\gamma_{4} + 2t^2\gamma_{3} - 2\gamma_{5}, \;  
     \rho_{6}   = \gamma_{3}^2 + 2c_{6} - 2t\gamma_{5} - 3t^2 \gamma_{4} + t^6, \\
     \rho_{8}  &= 3 \gamma_{4}^2 - 2\gamma_{3} \gamma_{5} + t(2c_{7} - 6\gamma_{3}\gamma_{4}) - 9t^2 c_{6} 
                  + 12t^3 \gamma_{5} + 15t^4 \gamma_{4} - 6t^5 \gamma_{3} - t^8, \\
     \rho_{9}  &= 2c_{6}\gamma_{3} + t^2 c_{7}  - 3t^3 c_{6} - 2\gamma_{9}, \; 
     \rho_{10}  = \gamma_{5}^2 - 2c_{7}\gamma_{3} + 3t^3c_{7},  \\
     \rho_{12} & = 3c_{6}^2 - 2 \gamma_{4}^3 - 2c_{7} \gamma_{5} + 2\gamma_{3}\gamma_{4}\gamma_{5}  
                  + t(4c_{7}\gamma_{4} - 2c_{6}\gamma_{5} + 6\gamma_{3} \gamma_{4}^2)  
                  + t^2(-3c_{7}\gamma_{3} + 3c_{6}\gamma_{4})   \\
               &  + t^3(-12\gamma_{4}\gamma_{5} + 5c_{6}\gamma_{3}) + t^4(-2\gamma_{3}\gamma_{5} - 15 \gamma_{4}^2)  
                  - 10t^6c_{6} + 12t^7\gamma_{5}+ 19t^8 \gamma_{4} - 6t^9 \gamma_{3} - 2t^{12}, \\
    \rho_{14}  & = c_{7}^2 + 6c_{7}\gamma_{3}\gamma_{4} - 2c_{6}\gamma_{3}\gamma_{5} - t^2c_{7}\gamma_{5}
                  + t^3(-9c_{7}\gamma_{4} + 3c_{6}\gamma_{5})  \medskip  
                  - 6t^4 c_{7}\gamma_{3}   + 9t^7 c_{7},  \medskip \\
   \rho_{18}   & =  - \gamma_{9}^2 + 2c_{6}c_{7}\gamma_{5} + 6c_{7}\gamma_{3} \gamma_{4}^2 - 2 c_{7}^2 \gamma_{4} 
                   -2c_{6}\gamma_{3}\gamma_{4}\gamma_{5} + 2c_{6}\gamma_{3}\gamma_{9} 
                  + t(-6 c_{7}^2\gamma_{3} + 24c_{6}c_{7}\gamma_{4})  \\
                & + t^2(-25c_{7}\gamma_{4}\gamma_{5} + c_{7}\gamma_{9} - 18c_{6}c_{7}\gamma_{3}) 
                  + t^3(-45c_{7} \gamma_{4}^2 + 20c_{7}\gamma_{3}\gamma_{5} + 3c_{6}\gamma_{4}\gamma_{5} 
                   - 3c_{6}\gamma_{9})  \medskip   \\
               & + t^4(11c_{7}^2 + 2c_{6}\gamma_{3}\gamma_{5}  + 48c_{7}\gamma_{3}\gamma_{4}) 
                + 51t^5 c_{6}c_{7}   - 53t^6 c_{7}\gamma_{5}  
                + t^7 (-69c_{7}\gamma_{4} - 3c_{6}\gamma_{5})  \\
               &+ 16t^8 c_{7}\gamma_{3} + 15t^{11} c_{7}.  
    \end{align*}
\end{thm}
\begin{rem}  \rm 
Similarly to  the case of $E_{6}$, 
this description is slightly different from the one in \cite{Nak01}.
The second author  \cite{Nak01}  make   use of the integral cohomology ring of the homogeneous space 
$E_7/T^1 \! \cdot \! \mathit{Spin}(12)$ to determine the higher relations $\rho_{12}$, $\rho_{14}$, and $\rho_{18}$.
Here we give  a description directly related to the ring of invariants of the Weyl group 
$W(E_{7})$ (see also \S \ref{Borel_presentation}).  
\end{rem}

In a similar fashion, we can give a Borel presentation  of $H^*(E_{8}/T;\Z)$.
\begin{thm}[\cite{Nak09.2}, Theorem 3.4]    \label{thm:E_8/T}   
  The integral cohomology ring of $E_{8}/T$ is 
     \begin{align*} 
             H^{*}(E_{8}/T;\Z) =\  &  \Z [t_{1}, \ldots, t_{8}, t, 
                                            \gamma_{3}, \gamma_{4}, \gamma_{5}, \gamma_{6}, \gamma_{9},
                                            \gamma_{10}, \gamma_{15} ]  \\
                                    & /(\rho_{1}, \rho_{2}, \rho_{3}, \rho_{4}, \rho_{5}, \rho_{6},  \rho_{8}, \rho_{9}, 
                                        \rho_{10}, \rho_{12}, \rho_{14}, \rho_{15}, \rho_{18}, \rho_{20}, \rho_{24}, \rho_{30}),   
    \end{align*}
    where
\begin{align*}
     \rho_{1} &= c_{1}-3t, \;  
     \rho_{2}  = c_{2}-4t^{2}, \;   
     \rho_{3}  = c_{3}-2\gamma_{3}, \; 
     \rho_{4}  = c_{4}+2t^{4}-3\gamma_{4}, \\
     \rho_{5} &= c_{5}-3t\gamma_{4} + 2t^{2}\gamma_{3}-2\gamma_{5},  \; 
     \rho_{6}  = c_{6}-2\gamma_{3}^{2}-t\gamma_{5}+t^{2}\gamma_{4}-2t^{6}-5\gamma_{6},   \\
     \rho_{8} &= -3c_{8}+3\gamma_{4}^{2}-2\gamma_{3}\gamma_{5}+ t(2c_{7}-6\gamma_{3}\gamma_{4}) 
                 + t^2(2\gamma_{3}^{2}-5\gamma_{6}) 
                 +3t^{3}\gamma_{5}+4t^{4}\gamma_{4}-6t^{5}\gamma_{3}+t^{8},  \\
     \rho_{9}  &= 2c_{6}\gamma_{3}+tc_{8}+t^{2}c_{7}-3t^{3}c_{6}-2\gamma_{9}, \; 
     \rho_{10}  = \gamma_{5}^{2}-2c_{7}\gamma_{3}  -t^{2}c_{8}  +3t^{3}c_{7}  -3\gamma_{10}, \\
     \rho_{12} &= 15 \gamma_{6}^2 + 2 \gamma_{3} \gamma_{4} \gamma_{5} -2 c_{7} \gamma_{5} + 2 \gamma_{3}^4 
                  + 10 \gamma_{3}^{2} \gamma_{6}  - 3 c_{8} \gamma_{4} - 2 \gamma_{4}^3 
                  + t (c_{8} \gamma_{3} -2 \gamma_{3}^2 \gamma_{5} + 4 c_7 \gamma_4 + 6 \gamma_{3} \gamma_{4}^2)   \\  
               &+ t^2 (3 \gamma_{10} - 25 \gamma_4 \gamma_6 - c_7 \gamma_3 -16  \gamma_{3}^2 \gamma_4)
                + t^3 (25 \gamma_3 \gamma_6 - 3 \gamma_4 \gamma_5 + 10  \gamma_{3}^3) 
                 + t^4 (3c_8 + 3 \gamma_3 \gamma_5 + 5 \gamma_{4}^2)   \\
               & + t^5 (-3 c_7 - 5 \gamma_3 \gamma_4)  + 4 t^6  \gamma_{3}^2 - 7 t^8 \gamma_4 + 4 t^9 \gamma_3, \\
     \rho_{14} &= c_{7}^2 - 3 c_8 \gamma_6 + 6 \gamma_4 \gamma_{10} - 4c_8 \gamma_{3}^2 + 6 c_{7} \gamma_3 \gamma_4  
                  - 6  \gamma_{3}^2 \gamma_{4}^2   - 12 \gamma_{4}^2 \gamma_6 - 2 \gamma_3 \gamma_5 \gamma_6 \\
               &  + t(24 \gamma_3 \gamma_4 \gamma_6 - 8c_{7}  \gamma_{3}^2  - 8 c_7 \gamma_6 + 4 c_{8} \gamma_5 
                  - 6 \gamma_{3} \gamma_{10} + 12  \gamma_{3}^3 \gamma_4)  \\
               &  + t^2 (-2 \gamma_3 \gamma_4 \gamma_5 + 6  \gamma_{4}^3 + 2  \gamma_{3}^2 \gamma_6 + 20  \gamma_{6}^2 
                          - 4  \gamma_{3}^4   - c_7 \gamma_5)   
                  + t^3 (-12 \gamma_3  \gamma_{4}^2 + 8c_8 \gamma_3 - 5 c_7 \gamma_4 + 3\gamma_5 \gamma_6) \\
               & + t^4 (3 \gamma_{10} - 26 \gamma_4 \gamma_6 + 6c_7 \gamma_3 - 4 \gamma_{3}^2 \gamma_4) 
                 + t^5 (24 \gamma_3 \gamma_6 + 3 \gamma_4 \gamma_5 + 12  \gamma_{3}^3) 
                  + t^6 (-6 c_8 + 2 \gamma_{4}^2)  \\
              & - 4 t^7 c_7 + t^8 (6 \gamma_6 - 6  \gamma_{3}^2) 
                - 6 t^{10} \gamma_4 + 12 t^{11} \gamma_3 - 2t^{14}, \\
     \rho_{15} &= (c_8 - t^2 c_6 + 2 t^3 \gamma_5 + 3 t^4 \gamma_4 - t^8)(c_7 - 3tc_6) 
                - 2(\gamma_{3}^2 + c_6)(\gamma_9 - c_6 \gamma_3)  - 2\gamma_{15}, 
\end{align*} 
\begin{align*} 
    \rho_{18}  &=  \gamma_{9}^2 - 9c_{8}\gamma_{10} - 6\gamma_{4}^2\gamma_{10} - 4\gamma_{3}^3\gamma_{9} 
                    - 10  \gamma_{3} \gamma_{6}\gamma_{9} + 2\gamma_{3}\gamma_{5}\gamma_{10}  
                   - 2\gamma_{3}\gamma_{4}\gamma_{5}\gamma_{6} 
                   - 6c_{7}\gamma_{3}\gamma_{4}^2 + 3c_{8}\gamma_{4}\gamma_{6}  \\
               & + c_{8}\gamma_{3}^2\gamma_{4}
                   + 6\gamma_{3}^2\gamma_{4}^3   + 12\gamma_{4}^3 \gamma_{6}+ 2c_{7}^2\gamma_{4} 
                   + 2c_{7}\gamma_{3}^2\gamma_{5}   - 2\gamma_{3}^3\gamma_{4}\gamma_{5}+ 2c_{7}\gamma_{5}\gamma_{6} 
                   + 4\gamma_{3}^6  
                  - 10 \gamma_{6}^3  + 18\gamma_{3}^4 \gamma_{6}  \\
               & + 15 \gamma_{3}^2 \gamma_{6}^2  - 9c_{7}c_{8}\gamma_{3} \\
              & + t(-2\gamma_{3}\gamma_{5}\gamma_{9}  - 24c_{7}\gamma_{4} \gamma_{6} + 8c_{8}\gamma_{4}\gamma_{5}
                + 4c_{7}\gamma_{3}^2\gamma_{4} 
                + 4c_{7}\gamma_{10}   - c_{8}\gamma_{9} + 2c_{7}^2 \gamma_{3}  + 4c_{8}\gamma_{3} \gamma_{6} 
                + 12 \gamma_{3} \gamma_{4} \gamma_{10} \\
              &  - 36\gamma_{3}\gamma_{4}^2 \gamma_{6}   + 12 \gamma_{3}^2\gamma_{5} \gamma_{6} 
                + c_{8}\gamma_{3}^3 + 6\gamma_{3}^4\gamma_{5} - 18\gamma_{3}^3 \gamma_{4}^2)  \\
             &  + t^{2} (24\gamma_{3}^4\gamma_{4} - 2c_{8}^2 - c_{7}\gamma_{9} - 11\gamma_{3}^2\gamma_{10} 
                   + 2\gamma_{3}\gamma_{4}\gamma_{9}   - 2c_{8}\gamma_{3}\gamma_{5}   + 16c_{7}\gamma_{3} \gamma_{6} 
                   - 3c_{7}\gamma_{4}\gamma_{5}  
                + 75\gamma_{4} \gamma_{6}^2   \\  
             & - 6\gamma_{4}^4  - 9c_{8}\gamma_{4}^2      + 81 \gamma_{3}^2\gamma_{4} \gamma_{6} - 13 \gamma_{6}\gamma_{10} 
                   + 4\gamma_{3}\gamma_{4}^2 \gamma_{5}  - c_{7}\gamma_{3}^3)  \\
               & + t^3(-3\gamma_{5}\gamma_{10} - 150\gamma_{3} \gamma_{6}^2 - 135\gamma_{3}^3 \gamma_{6}
                 + 6\gamma_{3}^2 \gamma_{9}  
                 - 2 c_{7}\gamma_{3}\gamma_{5}   + 21c_{7}\gamma_{4}^2  + 15c_{7}c_{8}   + 3\gamma_{4}\gamma_{5} \gamma_{6} 
                      - 3\gamma_{3}^2 \gamma_{4} \gamma_{5}  \\
             &  + 18\gamma_{3} \gamma_{4}^3   + 15 \gamma_{6}\gamma_{9} 
                       + 14c_{8} \gamma_{3} \gamma_{4} - 30\gamma_{3}^5)   \\
               & + t^4(-13c_{8} \gamma_{6} + 2\gamma_{4}\gamma_{10} - 5c_{7}^2 - 33\gamma_{3}^2 \gamma_{4}^2 
                 + 3\gamma_{5} \gamma_{9}   -28 \gamma_{3}\gamma_{5} \gamma_{6} - 45\gamma_{4}^2 \gamma_{6}  
                  - 41c_{7}\gamma_{3}\gamma_{4} -13 \gamma_{3}^3 \gamma_{5}  \\
               & - 9c_{8} \gamma_{3}^2) \\
               & + t^5(3c_{7} \gamma_{6} - 6\gamma_{4}^2\gamma_{5} + 23c_{7}\gamma_{3}^2  
                 + 105\gamma_{3}\gamma_{4} \gamma_{6} - 6c_{8} \gamma_{5}  - 3\gamma_{4}\gamma_{9} 
                 + 45 \gamma_{3}^3 \gamma_{4})   \\
               & + t^6(11\gamma_{4}^3 - 4\gamma_{3}\gamma_{9} + 4c_{7}\gamma_{5} + 9\gamma_{3}\gamma_{4}\gamma_{5} 
                 + 12\gamma_{3}^4    + 66\gamma_{3}^2 \gamma_{6}  + 75\gamma_{6}^2 + 2c_{8}\gamma_{4})  \\
               & + t^7(-33\gamma_{3}\gamma_{4}^2 + 12 \gamma_{3}^2 \gamma_{5}  + 15\gamma_{5} \gamma_{6})  
                 + t^8(-4\gamma_{10} + 21\gamma_{3}^2 \gamma_{4} - 5c_{7}\gamma_{3} -3\gamma_{4} \gamma_{6}) \\
               & + t^9(6\gamma_{9}  - 42\gamma_{3}^3 - 99\gamma_{3} \gamma_{6}) 
                 + t^{10} (-4c_{8} - 6\gamma_{4}^2 -13 \gamma_{3}\gamma_{5})
                 + t^{11}(3c_{7} + 27\gamma_{3}\gamma_{4}) 
                 +t^{12} (60 \gamma_{6} + 18\gamma_{3}^2)  \\
               & + 6t^{13} \gamma_{5} - 9t^{14} \gamma_{4} 
                 - 12t^{15} \gamma_{3} +10t^{18}, \\
         \rho_{20}   &= 9u^{20} + 45u^{14}v + 12u^{10}w + 60u^8v^2 + 30u^4vw  + 10u^2v^3 + 3w^2, \\
         \rho_{24}   &=  11u^{24} + 60u^{18}v +  21u^{14}w + 105u^{12}v^2 + 60u^8vw 
                        + 60u^{6}v^3 + 9u^4 w^2  + 30 u^2 v^2 w + 5v^4,   \\
         \rho_{30}   &= -9x^2 - 12u^{9}vx  -6u^{5}wx  + 9u^{14}vw -10u^{12}v^3  -3u^{10}w^2  + 30u^8 v^2w  
                        -35u^{6}v^{4} + 6u^{4}vw^{2}  \\
                     & - 10u^2v^3w  -4v^5 -2w^{3}, 
     \end{align*}
and 
  \begin{align*}
      u & = t_{8}, \\
      v & =  2\gamma_{6}+ \gamma_{3}^{2}-u\gamma_{5}+\gamma_{4}(-t^{2}+u^{2})
            -u^{3}\gamma_{3} + t^{6} - t^{4}u^{2} + t^{3}u^{3} + t^{2}u^{4} - tu^{5}, \\
      w & = \gamma_{10} + u\gamma_{9} - u^{3}c_{7}  -u\gamma_{4}\gamma_{5} 
            +2u^{2} \gamma_{4}^{2} -2u^{2}\gamma_{3}\gamma_{5} 
            + \gamma_{3}\gamma_{4}(-6tu^{2} + 2u^{3}) +  \gamma_{3}^{2}(2t^{2}u^{2} + 2tu^{3} - 2u^{4})  \\
        &  + \gamma_{6}(-5t^{2}u^{2} + 5tu^{3})
         + \gamma_{5}(t^{4}u + 3t^{3}u^{2} + t^{2}u^{3})  
         +\gamma_{4}(6t^{4}u^{2} - 3t^{3}u^{3} -2t^{2}u^{4} -tu^{5} + u^{6}) \medskip \\ 
        & +\gamma_{3}(-6t^{5}u^{2} - 2t^{4}u^{3} + 4t^{3}u^{4} + 6t^{2}u^{5} - 4tu^{6} + u^{7}) 
          + 4t^{7}u^{3} -6t^{5}u^{5} + 2t^{4}u^{6}  + t^{3}u^{7} - t^{2}u^{8}, \\
     x &  =  \gamma_{15}-20 \gamma_{3} \gamma_{6}^2 +3  \gamma_{3}^{2}\gamma_{9}  -23 \gamma_{3}^3 \gamma_{6} 
            -6 \gamma_{3}^{5}  +4 \gamma_{6}\gamma_{9} +3 u \gamma_{4}\gamma_{10} - u \gamma_{5}\gamma_{9}  
           -3 u \gamma_{3}^{2} \gamma_{4}^{2} +3 uc_{7}\gamma_{3}\gamma_{4}  \\
         & -6 u \gamma_{4}^2\gamma_{6}  + (-3t + 2 u ) \gamma_{3}^{3} \gamma_{5}   
           +(-4 t  + 4u) \gamma_{3}\gamma_{5}\gamma_{6}  
           +(-t^2 - u^2) \gamma_{4}\gamma_{9}  +(t^2 + tu -u^2 )c_{7}   \gamma_{3}^{2}  \\
           & +(9t^2+ 12 tu + 5u^2) \gamma_{3}\gamma_{4}\gamma_{6} 
    +  (5t^2 +6tu + 2 u^2)  \gamma_{3}^{3} \gamma_{4}  + (3 t^2 +4tu + u^2)  c_{7} \gamma_{6}        \\
  & +(-6t^3 -2t^2u - 6tu^2 + 5u^3) \gamma_{3}^{4}    - u^3 \gamma_{3}\gamma_{9} 
    + (3t^2 u + u^3 ) \gamma_{4}^{3}   + (2 t^2u + 3tu^2) c_{7} \gamma_{5}  \\
  & +(-45t^3 + 10t^2 u -40 tu^2) \gamma_{6}^{2} 
    + (t^3 -2t^2 u + tu^2 - u^3 )\gamma_{3} \gamma_{4} \gamma_{5}  
    +(-33 t^3 + t^2u - 31tu^2 + 13u^3) \gamma_{3}^{2} \gamma_{6} \\
  & +( -2 t^4 - 4 t^3u -3t u^3 + 3u^4) c_{7} \gamma_{4}   
    +(-9 t^4 -6t^3u -18t^2u^2 + 5tu^3 -3u^4)\gamma_{5} \gamma_{6}  \\
  & +(-3t^4 -3t^3 u -7t^2 u^2 + 5tu^3 -4u^4) \gamma_{3}^{2}\gamma_{5} 
    +(-t^4 -6 t^3u  -t^2u^2 -3tu^3)\gamma_{3} \gamma_{4}^2    
\end{align*} 
\begin{align*} 
  & + (-3t^4 u -6 t^3 u^2 + 3t^2 u^3  +15 t u^4  )\gamma_{10} 
    + (-3t^4 u +  t^3 u^2 + 5 t^2 u^3 + 10tu^4  -u^5 )c_{7} \gamma_{3}   \\
  &  +(15 t^5 -2t^4u + 3t^3 u^2 + 14t^2 u^3 -16tu^4 + 3u^5) \gamma_{3}^2 \gamma_{4}  \\
  &   + (39t^5 - 13t^4u + 8t^3 u^2 + 35 t^2 u^3 -31tu^4 - 3u^5)\gamma_{4} \gamma_{6}    
     +(t^6 -t^4u^2 -{t}^{3}{u}^{3} - t^2 u^4 - tu^5 -u^6)\gamma_{9}    \\
  & +(- 13t^6 + 12t^5 u +5t^4 u^2 -56t^3 u^3 + 8t^2 u^4 +21tu^5 + 2u^6) \gamma_{3} \gamma_{6}       \\
  & + ( 6 t^6 + 3t^5 u +2t^4 u^2 +  7 t^3u^3  + t^2 u^4 -8tu^5  + 3u^6)\gamma_{4}  \gamma_{5}    \\
  & +(-8 t^6 + 6t^5u + 2t^4u^2 -22t^3 u^3 + 6t^2u^4 +8tu^5   -2u^6) \gamma_{3}^{3}   \\
  & + (-6t^7 + t^6 u -7t^4 u^3  + 5 t^3 u^4  + 3t^2 u^5 + 3tu^6  - 63 u^7)  \gamma_{4}^{2}   \\
  & + (-t^7 +2t^6 u + t^5 u^2 -11 t^4 u^3 + 6t^3 u^4 + 5t^2 u^5 +6tu^6   + 39u^7)\gamma_{3}\gamma_{5}   \\
  & + (2t^8 + 6t^7 u +3t^6u^2 -4t^5 u^3  - 15 t^4 u^4  + 6t^3 u^5   + 3t^2 u^6  -40tu^7 + 59 u^8)c_{7}   \\
  & + (3t^8 +  t^6 u^2 +11t^5 u^3 + 14t^4 u^4 -20t^3 u^5 -4 t^2 u^6  +118tu^7 + 3u^8 )\gamma_{3} \gamma_{4}    \\  
  & + (- 48 t^9 + 3 t^8u  -41 t^7 u^2  + 18 t^6 u^3 + 16 t^5 u^4 
    -13t^4 u^5  - 67t^3 u^6 + 125t^2 u^7 - 15tu^8  -291u^9)\gamma_{6}    \\
  & +(-18t^9 -3t^8 u -16t^7 u^2 + 10t^6u^3 -4t^5u^4 -8t^4 u^5  
    -16t^3 u^6  -23t^2 u^7 -10tu^8  -115u^9)  \gamma_{3}^{2} \\
  &  +(-6t^{10} -3t^9 u -9t^8 u^2 + 5t^7 u^3 -5 t^6 u^4  -14t^4 u^6  
     -52t^3 u^7 + 6t^2 u^8 -60tu^9 +117 u^{10})\gamma_{5}   \\
  & + (18t^{11} -3t^{10} u+ 5t^9 u^2 + 11t^8 u^3 - 28t^7 u^4 + 8t^6 u^5  
    + 20t^5 u^6  -64t^4 u^7 -15t^3 u^8  + 54t^2 u^9  \\
  & + 178t u^{10}   - 177u^{11}) \gamma_{4} \\
  & +(-2t^{12} +6t^{11}u + 2t^{10}u^2 -20 t^9 u^3 +11t^8 u^4  + 22t^7 u^5 -8t^6 u^6 + 83t^5 u^7  + 15t^4 u^8 + 5t^3 u^9 \\
  & -116 t^2 u^{10}  
    + tu^{11} + 117 u^{12}) \gamma_{3}        \\
 &  -12t^{15} - t^{14}u -10t^{13}u^2 + 6 t^{12} u^3 + 7t^{11}u^4   -13t^{10}u^5  -31t^9 u^6  +9 t^8 u^7 -t^7 u^8 -118t^6 u^9    
    -18t^5 u^{10}  \\
 &+ 131t^4 u^{11}  -6t^3 u^{12}   - 233t^2u^{13}  + 175tu^{14} -58 u^{15}.        
  \end{align*} 
\end{thm}  
\begin{rem}  \rm 
  In order to determine the higher relations $\rho_{20}$, $\rho_{24}$, and $\rho_{30}$, the second 
  author \cite{Nak09.2} make use of the integral cohomology ring of the homogeneous space $E_{8}/T^{1} \! \cdot \! E_{7}$
  (for the integral cohomology ring of $E_{8}/T^{1} \! \cdot \! E_{7}$, see \cite{Nak09}). 
   
\end{rem} 

Since the elements $t_{i} \; (1 \leq i \leq l)$ in  $H^{2}(E_{l}/T;\Z)$ belong to an orbit of 
$W_{P_{2}}$, the elementary symmetric polynomials $c_{i} \; (1 \leq i \leq l)$  are $W_{P_{2}}$-invariants, 
and hence they are elements of $H^{*}(E_{l}/P_{2};\Z)$.  The element $t$ in $H^{2}(E_{l}/T;\Z)$ 
is also an element of $H^{2}(E_{l}/P_{2};\Z)$.   Therefore,  in any cases of $E_{l} \ (l=6,7,8)$,  
one can see that the higher degree generators $\gamma_{i}$  lie in $H^*(E_{l}/P_2;\Z) \subset H^*(E_{l}/T;\Z)$.

\subsection{}   
In the previous subsection, we reviewed  the Borel presentation of the integral cohomology rings of 
$E_{l}/T$ by   generators and relations. In order to  express those generators 
in terms  of Schubert classes, we will use   the divided difference operators  introduced  independently  
by Bernstein-Gelfand-Gelfand  \cite{BGG73} and  Demazure  \cite{Dem73}:  
\begin{df}[Bernstein-Gelfand-Gelfand \cite{BGG73}, Demazure  \cite{Dem73}]  
$(1)$  For each root $\alpha \in \Delta$, the operator of degree $-2$
   \begin{equation*} 
        \Delta_{\alpha} :  H^*(BT;\Z) \longrightarrow H^*(BT;\Z) 
  \end{equation*} 
is defined as   
  \begin{equation*} 
        \Delta_{\alpha}(u) = \dfrac{u - s_{\alpha}(u)}{\alpha}  \quad \text{for} \quad 
         u  \in   H^*(BT;\Z).   
  \end{equation*} 
$(2)$ For $w \in W$,  the operator $\Delta_w$  is defined as  the composite 
    \begin{equation*} 
        \Delta_{w} = \Delta_{\alpha_{1}} \circ \Delta_{\alpha_{2}} \circ \cdots \circ \Delta_{\alpha_{k}},  
  \end{equation*} 
where $w = s_{\alpha_{1}} s_{\alpha_{2}} \cdots s_{\alpha_{k}} \ (\alpha_{i} \in \Pi)$ is any reduced 
decomposition of $w$. 
\end{df} 
 One can show that the definition is well defined, i.e., 
independent of the choice of the  reduced decomposition of $w$. 
The following theorem gives a correspondence between a cycle represented by a polynomial
in the Borel presentation and a sum of Schubert classes via the characteristic homomorphism (\ref{eqn:ch.hom2}).
\begin{thm}[Bernstein-Gelfand-Gelfand \cite{BGG73}, Demazure \cite{Dem73}] \label{thm:BGG} 
  For a  homogeneous polynomial $f \in H^{2k}(BT;\Z)$,  we have
  \begin{equation*} 
         c(f) = \sum_{w \in W, \ l(w) = k} \Delta_{w}(f) \,  Z_{w}.   
  \end{equation*}  
$($Note that $\Delta_{w}(f) \in H^{0}(BT;\Z) \cong \Z$.$)$  \;  In particular,  we have 
      $c(\omega_i) = Z_{s_i}$.
\end{thm}

\subsection{}   \label{subsec:BGG}  
Now we combine the result of the previous subsections to have the explicit relation
between the ring generators of $H^*(E_{l}/T;\Z)$
given in \S \ref{subsec:E_l/T} and the Schubert basis  $\{Z_{w} \}_{w \in W(E_{l})}$. 
For simplicity,  a reduced decomposition
$w = s_{i_{1}}s_{i_{2}}\cdots s_{i_{k}}$ will be  abbreviated to   $[i_{1} i_{2} \cdots i_{k}]$, 
and we denote  by  $Z_{i_{1}i_{2} \cdots i_{k}}$ the Schubert class corresponding to $w$, 
although the reduced decomposition of a Weyl group element may not  be unique.

   First we deal with the case of $E_{6}$. Since $c(\omega_{i}) = Z_{i}$, it follows immediately  from 
  (\ref{eqn:t_i}) that 
     \begin{equation*} 
           \begin{array}{llll} 
             t_{1} & \hspace{-0.2cm} = -Z_{1} + Z_{2}, \;  t_{2} = Z_{1} + Z_{2} - Z_{3}, \;  
             t_{3} = Z_{2} + Z_{3} - Z_{4}, \medskip \\
             t_{4} & \hspace{-0.2cm} = Z_{4} - Z_{5}, \;  t_{5}  = Z_{5} - Z_{6}, \;  t_{6} = Z_{6}, \;   
             t  = Z_{2}. \medskip \\
           \end{array} 
     \end{equation*}   
   For the higher degree generators,  we would like to have an expansion of $\gamma_i \; (i = 3, 4)$
   by the sum of Schubert classes.
   Note that they lie   in  $H^*(E_6/P_2;\Z)  \subset H^*(E_{6}/T;\Z)$, which has the Schubert basis indexed by 
   $W^{P_{2}}$.
  Since the length three  and four  elements of $W^{P_2}$ are  $\{[342], [542] \}$, 
    $\{ [1342], [3542], [6542] \}$ respectively, the elements
    $\gamma_{i} \; (i = 3,4)$ should be written as $\Z$-linear combinations  
\begin{align*} 
       \gamma_{3} & = a_{342} Z_{342}+  a_{542} Z_{542},  \\
       \gamma_{4} & = a_{1342} Z_{1342} + a_{3542} Z_{3542} + a_{6542} Z_{6542}. 
\end{align*} 
The coefficients $a_{342}, a_{542}, a_{1342}, a_{3542}, a_{6542}  \in \Z$ can be determined as follows.
 By Theorem \ref{thm:E_6/T}, we have 
      \begin{equation} \label{eqn:gamma_i(E_6)} 
              2\gamma_{3}   =  c_{3}, \quad 3\gamma_{4} =  c_{4} + 2t^4. 
      \end{equation} 
Therefore  $2\gamma_{3}$ and $3\gamma_{4}$ are contained in the image of $c$.
Define the polynomials in  $H^*(BT;\Z)$ by 
         \begin{equation*}
             \delta_{3} = c_{3}, \quad  \delta_{4}  = c_{4} + 2t^4, 
          \end{equation*} 
  so that they are mapped  under $c$ to $2\gamma_{3}$ and $3\gamma_{4}$ in $H^*(E_{6}/T;\Z)$ respectively.
Applying  Theorem \ref{thm:BGG} to $\delta_{i} \; (i = 3, 4)$, we have
    \begin{equation} \label{eqn:gamma_i(E_6).2}
           \begin{array}{rll} 
                  2\gamma_{3}   &\hspace{-0.2cm}  =  2Z_{342} + 4Z_{542} =  2(Z_{342} + 2Z_{542}),   \medskip  \\
                  3\gamma_{4}   &\hspace{-0.2cm}  =  3Z_{1342} + 6Z_{3542} + 6Z_{6542} 
                                                  =  3(Z_{1342} + 2Z_{3542} + 2Z_{6542}). \medskip 
            \end{array} 
   \end{equation} 
Since $H^*(E_{6}/T;\Z)$ is torsion free, we obtain    
       \begin{equation}  \label{eqn:gamma_i(E_6).3}
           \gamma_{3} = Z_{342}  + 2Z_{542}, \quad \gamma_{4} = Z_{1342} + 2Z_{3542} + 2Z_{6542}. 
       \end{equation} 
Conversely we can express each Schubert class as a polynomial in $t$, $\gamma_{3}$ and $\gamma_{4}$; Consider  the following polynomials in   $H^*(BT;\Q)$: 
   \begin{equation*} 
      \begin{array}{llll} 
          f_{342}  &\hspace{-0.2cm} =  -\dfrac{1}{2} \delta_{3}  + 2t^3, \;   
          f_{542}  = \dfrac{1}{2} \delta_{3} - t^3, \medskip  \\
          f_{1342} &\hspace{-0.2cm} = \dfrac{1}{3}\delta_{4} - t\delta_{3} + 2t^4, \;  
          f_{3542} = -\dfrac{1}{3}\delta_{4} + \dfrac{1}{2}t \delta_{3}, \; 
          f_{6542}   = \dfrac{1}{3}\delta_{4} - t^4. \medskip 
     \end{array}  
  \end{equation*}           
Applying Theorem \ref{thm:BGG}  to the above polynomials, we see that these are indeed the polynomial 
 representatives of the  Schubert classes, 
 i.e., $c(f_{342}) = Z_{342}$ and so on.   Therefore, in the ring $H^*(E_{6}/T;\Z)$, we obtain
  \begin{equation*} 
    \begin{array}{lll} 
       Z_{342}  &\hspace{-0.2cm} = -\gamma_{3} + 2t^3, \;   Z_{542} = \gamma_{3} - t^3, \medskip \\
       Z_{1342} &\hspace{-0.2cm} = \gamma_{4} - 2t\gamma_{3} + 2t^4, \;  Z_{3542} = -\gamma_{4} + t\gamma_{3}, \;  
       Z_{6542} = \gamma_{4} - t^4. \medskip 
    \end{array} 
   \end{equation*}    
Thus, for example, we can choose $Z_{542}, Z_{6542}$ as ring generators in place of $\gamma_{3}, \gamma_{4}$. 
Consequently,   we obtain the following result. 
\begin{prop}   \label{prop:E_6/T}
    In   $H^*(E_{6}/T;\Z)$, the relation between the ring generators $\{ t_{1}, \ldots, t_{6}, t,  
    \gamma_{3},  \gamma_{4} \}$  in Theorem $\ref{thm:E_6/T}$  and the Schubert classes is given by 
        \begin{equation*} 
            \begin{array}{llll} 
                t_{1} &\hspace{-0.2cm} =  -Z_{1} + Z_{2}, \;   t_{2} =  Z_{1} + Z_{2} - Z_{3}, \; 
                t_{3}  =  Z_{2} + Z_{3} - Z_{4}, \medskip   \\
                t_{4} &\hspace{-0.2cm}  =  Z_{4} - Z_{5},  \; 
                t_{5} = Z_{5} - Z_{6}, \;  t_{6} =  Z_{6}, \;  
                t      =  Z_{2}, \medskip \\
                \gamma_{3} &\hspace{-0.2cm} =  Z_{342}  + 2Z_{542}, \;  \gamma_{4} =  Z_{1342} + 2Z_{3542} + 2Z_{6542}.  \medskip  
           \end{array} 
        \end{equation*} 
Furthermore, we have
        \begin{equation*}  
               \begin{array}{lll}
                 Z_{342}  &\hspace{-0.2cm} =  -\gamma_{3}  + 2t^3, \;  Z_{542} =   \gamma_{3}  - t^3, \medskip  \\
                 Z_{1342} &\hspace{-0.2cm} = \gamma_{4} - 2t\gamma_{3} + 2t^4, \;   Z_{3542} = -\gamma_{4} + t\gamma_{3}, \; 
                 Z_{6542} = \gamma_{4} - t^4. \medskip 
               \end{array} 
        \end{equation*} 
In particular, the set   $\{ Z_{1}, Z_{2}, \ldots, Z_{6}, Z_{542}, Z_{6542} \}$ is  a minimal system of 
 ring generators of $H^*(E_{6}/T;\Z)$ that consists of Schubert classes. 
\end{prop}

  In a similar manner, we carry out the computation for  the case of $E_{7}$.  It follows  immediately  
  from (\ref{eqn:t_i}) that 
   \begin{equation*} 
       \begin{array}{llll} 
            t_{1} &\hspace{-0.2cm} = -Z_{1} + Z_{2}, \;  t_{2} = Z_{1} + Z_{2} - Z_{3}, \;  
            t_{3} = Z_{2} + Z_{3} - Z_{4}, \medskip \\
            t_{4} &\hspace{-0.2cm} = Z_{4} - Z_{5}, \; t_{5} = Z_{5} - Z_{6}, \; t_{6} = Z_{6} - Z_{7}, \; 
            t_{7} = Z_{7}, \;  t = Z_{2}.  \medskip 
        \end{array} 
     \end{equation*} 
By Theorem \ref{thm:E_7/T}, we have 
      \begin{equation} \label{eqn:gamma_i(E_7)}
          \begin{array}{llll} 
                2\gamma_{3}  &\hspace{-0.2cm} =  c_{3}, \medskip  \\
                3\gamma_{4}  &\hspace{-0.2cm} =  c_{4} + 2t^4,  \medskip  \\
                2\gamma_{5}  &\hspace{-0.2cm} =  c_{5} - 3t\gamma_{4} + 2t^2 \gamma_{3}  
                                              =  c_{5} - tc_{4} + t^2 c_{3} - 2t^5, \medskip  \\
                2\gamma_{9}  &\hspace{-0.2cm} =  2c_{6}\gamma_{3} + t^2 c_{7} - 3t^3 c_{6}  
                                              =  c_{3}c_{6} + t^2 c_{7} - 3t^3 c_{6}. 
           \end{array} 
      \end{equation} 
Therefore   $2\gamma_{3}$, $3\gamma_{4}$, $2\gamma_{5}$, and $2\gamma_{9}$  are contained in the image of $c$. 
Define the polynomials in $H^{*}(BT;\Z)$ by  
\[  
   \delta_{3} = c_{3}, \; \delta_{4}= c_{4} + 2t^{4},  \; \delta_{5}= c_{5} - tc_{4} + t^2 c_{3} - 2t^5, \; 
   \delta_{9}  = c_{3}c_{6} + t^2 c_{7} - 3t^3 c_{6}, 
\] 
so that they are mapped under $c$  to $2\gamma_{3}$, $3\gamma_{4}$, $2\gamma_{5}$, and $2\gamma_{9}$ in $H^*(E_{7}/T;\Z)$ 
respectively.   Applying  Theorem \ref{thm:BGG} to  $\delta_{i} \; (i = 3, 4, 5, 9)$, we have
    \begin{equation} \label{eqn:gamma_i(E_7).2}
           \begin{array}{rll} 
                 2\gamma_{3}   &\hspace{-0.2cm} =  2(Z_{342} + 2Z_{542}),   \medskip  \\
                 3\gamma_{4}   &\hspace{-0.2cm} =  3(Z_{1342} + 2Z_{3542} + 2Z_{6542}), \medskip   \\
                 2\gamma_{5}   &\hspace{-0.2cm} = 2Z_{76542}, \medskip \\
                 2\gamma_{9}   &\hspace{-0.2cm} = 2(2Z_{154376542} + Z_{654376542}), \medskip  
             \end{array}  
   \end{equation}  
Since $H^*(E_{7}/T;\Z)$ is torsion free, we obtain
   \begin{equation} 
           \begin{array}{rll} 
                 \gamma_{3}   &\hspace{-0.2cm} =  Z_{342} + 2Z_{542},   \medskip  \\
                 \gamma_{4}   &\hspace{-0.2cm} =  Z_{1342} + 2Z_{3542} + 2Z_{6542}, \medskip   \\
                 \gamma_{5}   &\hspace{-0.2cm} =  Z_{76542}, \medskip \\
                 \gamma_{9}   &\hspace{-0.2cm} =  2Z_{154376542} + Z_{654376542}, \medskip  
             \end{array}  
   \end{equation}   
Conversely we can express each Schubert class as a polynomial in $t$, $\gamma_{i} \; (i = 3, 4, 5, 9)$; 
Consider the following polynomials in  $H^*(BT;\Q)$: 
 \begin{equation*} 
  \begin{array}{llll} 
         & f_{342} =  -\dfrac{1}{2} \delta_{3}  + 2t^3, \;  f_{542}  = \dfrac{1}{2} \delta_{3}  - t^3, \medskip  \\
         & f_{1342} = \dfrac{1}{3} \delta_{4} - t\delta_{3} + 2t^4, \; 
           f_{3542} = -\dfrac{1}{3} \delta_{4} + \dfrac{1}{2} t \delta_{3}, \;  
           f_{6542} = \dfrac{1}{3} \delta_{4} - t^4, \medskip  \\
         & f_{76542} = \dfrac{1}{2} \delta_{5}, \medskip \\ 
         & f_{154376542} = \dfrac{1}{2} \delta_{9} - \dfrac{1}{6}\delta_{4}\delta_{5} + \dfrac{1}{2} t^4 \delta_{5},  \; 
           f_{654376542} = - \dfrac{1}{2} \delta_{9} + \dfrac{1}{3} \delta_{4}\delta_{5} - t^4 \delta_{5}.
                           \medskip   
   \end{array} 
  \end{equation*}
 Applying Theorem \ref{thm:BGG}  again to the above polynomials, we see that these are indeed the polynomial 
 representatives of the  Schubert classes, 
 i.e., $c(f_{342}) = Z_{342}$ and so on.   Therefore, in the ring $H^*(E_{7}/T;\Z)$, we obtain 
      \begin{equation*} 
         \begin{array}{lll} 
                & Z_{342} =  -\gamma_{3}  + 2t^3, \;  Z_{542} =   \gamma_{3}  - t^3, \medskip  \\
                & Z_{1342} = \gamma_{4} - 2t\gamma_{3} + 2t^4, \;  Z_{3542} = -\gamma_{4} + t\gamma_{3}, \; 
                  Z_{6542} = \gamma_{4} - t^4,   \medskip  \\
                & Z_{76542} =  \gamma_{5}, \medskip \\
                & Z_{154376542} = \gamma_{9} - \gamma_{4} \gamma_{5} + t^4 \gamma_{5}, \; 
                  Z_{654376542} = -\gamma_{9} + 2\gamma_{4} \gamma_{5} - 2t^4 \gamma_{5}. \medskip 
         \end{array} 
       \end{equation*} 
Consequently,   we obtain the following result. 
\begin{prop}   \label{prop:E_7/T}
    In $H^*(E_{7}/T;\Z)$, the relation between the ring generators $\{ t_{1}, \ldots, t_{7}, t, 
    \gamma_{3}, \gamma_{4}$,  $\gamma_{5}, \gamma_{9}  \}$  in Theorem $\ref{thm:E_7/T}$  and the Schubert classes is given by 
        \begin{equation*} 
            \begin{array}{llll}
                t_{1} & \hspace{-0.2cm} =  -Z_{1} + Z_{2}, \;  t_{2} =   Z_{1} + Z_{2} - Z_{3}, \; 
                t_{3} =  Z_{2} + Z_{3} - Z_{4}, \medskip \\
                t_{4} & \hspace{-0.2cm} =  Z_{4} - Z_{5}, \;  t_{5} =  Z_{5} - Z_{6}, \;  t_{6} =  Z_{6} - Z_{7}, \;  
                t_{7} =  Z_{7}, \;  t  =  Z_{2},  \medskip  \\
                         \gamma_{3} & \hspace{-0.2cm} =  Z_{342}  + 2Z_{542}, \;   
                         \gamma_{4} =  Z_{1342} + 2Z_{3542} + 2Z_{6542},   \medskip  \\
                         \gamma_{5} & \hspace{-0.2cm} =  Z_{76542}, \;  
                         \gamma_{9} =  2Z_{154376542} + Z_{654376542}. \medskip   
           \end{array} 
       \end{equation*} 
Furthermore, we have
        \begin{equation*}  
               \begin{array}{lll}
                 & Z_{342}  =  -\gamma_{3}  + 2t^3, \;  Z_{542} =   \gamma_{3}  - t^3, \medskip  \\
                 & Z_{1342} = \gamma_{4} - 2t\gamma_{3} + 2t^4, \;  Z_{3542} = -\gamma_{4} + t\gamma_{3}, \;   
                   Z_{6542} = \gamma_{4} - t^4, \medskip \\
                 & Z_{76542} =  \gamma_{5}, \medskip \\
                 & Z_{154376542} = \gamma_{9} - \gamma_{4} \gamma_{5} + t^4 \gamma_{5}, \; 
                   Z_{654376542} = -\gamma_{9} + 2\gamma_{4} \gamma_{5} - 2t^4 \gamma_{5}. \medskip 
               \end{array} 
        \end{equation*} 
In particular, the set   $\{ Z_{1}, Z_{2}, \ldots, Z_{7}, Z_{542}, Z_{6542}, Z_{76542}, 
Z_{654376542} \}$ is a minimal system of ring generators of $H^*(E_{7}/T;\Z)$ that consists of 
Schubert classes. 
\end{prop}

Finally, we carry out the computation for the case of $E_{8}$. 
It follows from (\ref{eqn:t_i}) that 
   \begin{equation*} 
       \begin{array}{llll} 
            t_{1} &\hspace{-0.2cm} = -Z_{1} + Z_{2}, \;  t_{2} = Z_{1} + Z_{2} - Z_{3}, \;  t_{3} = Z_{2} + Z_{3} - Z_{4}, \; 
            t_{4}  = Z_{4} - Z_{5},  \medskip \\
            t_{5} &\hspace{-0.2cm} = Z_{5} - Z_{6}, \; t_{6} = Z_{6} - Z_{7}, \; 
            t_{7} = Z_{7} - Z_{8}, \; t_{8} = Z_{8}, \;   t = Z_{2}.  \medskip 
        \end{array} 
     \end{equation*} 
By Theorem \ref{thm:E_8/T}, we have 
\begin{equation} \label{eqn:gamma_i(E_8)} 
  \begin{array}{lllll}  
    2\gamma_{3}    &= c_{3},  \medskip \\
    3\gamma_{4}    &= c_{4} + 2t^4, \medskip \\
    2\gamma_{5}    &= c_{5} - tc_{4} + t^2 c_{3} - 2t^5, \medskip \\
    30 \gamma_{6}  &= 6c_{6} - 3c_{3}^2 - 3tc_{5} + 5t^2 c_{4} - 3t^3 c_{3} - 2t^6, \medskip \\
    2\gamma_{9}    &= c_{3}c_{6} + tc_{8} + t^2 c_{7} - 3t^3 c_{6}, \medskip \\
    12\gamma_{10}  &=   (c_{5} - tc_{4}+ t^2c_{3}- 2t^5)^2 - 4c_{7}c_{3} - 4t^2 c_{8} + 12 t^3 c_{7}, \medskip \\
     8 \gamma_{15} &= 4(c_{8} - t^2 c_{6} + t^3 c_{5} + t^5 c_{3} - t^8)(c_{7} - 3tc_{6})  \\
                   &         - (c_{3}^2 + 4c_{6})(tc_{8} + t^2c_{7} -3t^3c_{6}).   
  \end{array} 
\end{equation} 
Therefore $2\gamma_{3}$, $3\gamma_{4}$, $2\gamma_{5}$, $30\gamma_{6}$, $2\gamma_{9}$, $12\gamma_{10}$, and $8\gamma_{15}$
are contained in the image of $c$.  Define the polynomials in $H^{*}(BT;\Q)$ by  
  \begin{align*}   
    \delta_{3}   & = c_{3},  \;  \delta_{4} = c_{4} + 2t^4, \; \delta_{5} = c_{5} - tc_{4} + t^2 c_{3} - 2t^5, \medskip \\
    \delta_{6}   & = \dfrac{1}{6}(6c_{6} - 3c_{3}^2 - 3tc_{5} + 5t^2 c_{4} - 3t^3 c_{3} - 2t^6), \medskip \\
    \delta_{9}   & = c_{3}c_{6} + tc_{8} + t^2 c_{7} - 3t^3 c_{6}, \medskip \\
    \delta_{10}  & =  \dfrac{1}{4} \{  (c_{5} - tc_{4}+ t^2c_{3}- 2t^5)^2 - 4c_{7}c_{3} - 4t^2 c_{8} + 12 t^3 c_{7} 
                                   \}, \medskip \\
    \delta_{15}  & = (c_{8} - t^2 c_{6} + t^3 c_{5} + t^5 c_{3} - t^8)(c_{7} - 3tc_{6})  
                     - \dfrac{1}{4}(c_{3}^2 + 4c_{6})(tc_{8} + t^2c_{7} -3t^3c_{6}).     
  \end{align*}  
By applying Theorem \ref{thm:BGG} to $\delta_{i} \; (i = 3, 4, 5, 6, 9, 10, 15)$,  we obtain the similar results: 
\begin{prop}   \label{prop:E_8/T}  
In $H^*(E_{8}/T;\Z)$, the relation between the ring gererators $\{ t_{1}, \ldots, t_{8}, t, \gamma_{3}, 
\gamma_{4}, \gamma_{5}$, $\gamma_{6},  \gamma_{9}, \gamma_{10}, \gamma_{15} \}$ in Theorem $\ref{thm:E_8/T}$
and the Schubert classes is  given by  

\begin{align*}  
            t_{1} & = -Z_{1} + Z_{2}, \;  t_{2} = Z_{1} + Z_{2} - Z_{3}, \;  t_{3} = Z_{2} + Z_{3} - Z_{4}, \; 
            t_{4}  = Z_{4} - Z_{5}, \\
            t_{5} & = Z_{5} - Z_{6}, \; t_{6} = Z_{6} - Z_{7}, \; 
            t_{7} = Z_{7} - Z_{8}, \; t_{8} = Z_{8}, \;   t = Z_{2}, \\
            \gamma_{3}  &=  Z_{342}   + 2Z_{542}, \\
      \gamma_{4}  &=  Z_{1342}  + 2Z_{3542}  + 2Z_{6542}, \\
      \gamma_{5}  &=  Z_{76542}, \\
      \gamma_{6}  & = -6Z_{136542} -5Z_{143542} -2Z_{243542} -5Z_{376542} -6Z_{436542} 
                      -Z_{876542}, \\
      \gamma_{9}  &=  4Z_{143876542}  + 2Z_{154376542}   + 4Z_{543876542}   +  Z_{654376542},   \\
      \gamma_{10}  &=  -Z_{15438765432},  \\
      \gamma_{15}  &=  58Z_{131426543876542} -17Z_{134231543876542} + 140Z_{134276543876542} \\
                   &+ 30Z_{135426543876542} + 127Z_{154276543876542} -22Z_{234231543876542} \\
                   &+ 87Z_{242316543876542} + 271Z_{243176543876542} -22Z_{245423143876542} \\
                   &+52Z_{254316543876542} +386Z_{314276543876542} +82Z_{315426543876542} \\
                   &+102Z_{342316543876542} -22Z_{345423143876542} +30Z_{354276543876542} \\
                   &+470Z_{423176543876542} -17Z_{458765423143542} +55Z_{465423143876542} \\
                   &+139Z_{542316543876542} +62Z_{543176543876542} +15Z_{654276543876542} \\
                   &+8Z_{658765423143542} +157Z_{765423143876542}.  
\end{align*} 
Furthermore, we have 
\begin{align*} 
   Z_{542}   &= \gamma_{3} - t^{3}, \\
   Z_{6542}  &= \gamma_{4} - t^{4}, \\
   Z_{76542} &= \gamma_{5}, \\
   Z_{136542} &= \gamma_{6} - t\gamma_{5} + t^{2} \gamma_{4}, \\ 
   Z_{154376542}  &  = \gamma_{9}  - 2\gamma_{3}^{3} - 4\gamma_{3}\gamma_{6}        
                     - \gamma_{4} \gamma_{5} + t(-6\gamma_{4}^{2} + 5c_{8} + 4 \gamma_{3}\gamma_{5}) 
                     + t^{2} (- 4c_{7} +14\gamma_{3} \gamma_{4})  \\
                  &  + t^{3} (-2\gamma_{3}^{2}+ 14\gamma_{6}) - 5 t^{4} \gamma_{5} 
                     - 10 t^{5} \gamma_{4} + 10t^{6}\gamma_{3}, \\
   Z_{1654376542} &= -\gamma_{10} + \gamma_{5}^{2}  -2\gamma_{3}^{2}\gamma_{4} -4\gamma_{4}\gamma_{6} 
                     + 2t^{2}\gamma_{4}^{2} +  t^{4}(2\gamma_{3}^{2} + 4\gamma_{6}) - 4t^{6}\gamma_{4} + 2t^{10}, \\
   Z_{134276543876542} &= \gamma_{15} + 16\gamma_{6}\gamma_{9} + \gamma_{5}\gamma_{10} + 6\gamma_{3}\gamma_{6}^{2} 
                           + 5\gamma_{3}^{2}\gamma_{9}  -\gamma_{3}^{3}\gamma_{6} + 4\gamma_{3}^{2}\gamma_{4}\gamma_{5}
                           + \gamma_{3}^{5} \\
                       &-12 \gamma_{3}\gamma_{5}c_{7} - 29\gamma_{3}\gamma_{4}c_{8}  \\
                       &  + t(- 167\gamma_{6}c_{8} - 6\gamma_{5}\gamma_{9} + 165\gamma_{4}^{2}\gamma_{6} 
                              - 96\gamma_{3}\gamma_{5}\gamma_{6}  + 32\gamma_{3}\gamma_{4} c_{7}  - 258\gamma_{3}^{2}c_{8}  \\
                       &  + 276\gamma_{3}^{2}\gamma_{4}^{2} -181 \gamma_{3}^{3}\gamma_{5} )  \\
                       &  + t^{2}(107\gamma_{6}c_{7} + 11\gamma_{5}c_{8}  + 93\gamma_{4}\gamma_{9} 
                          + 48\gamma_{3}\gamma_{10}  - 6\gamma_{4}^{2}\gamma_{5} - 945\gamma_{3}\gamma_{4}\gamma_{6} 
                          + 190\gamma_{3}^{2}c_{7} \\
                       & - 795\gamma_{3}^{3}\gamma_{4})  \\
                       &  +t^{3}(3\gamma_{6}^{2} - 31\gamma_{5}c_{7}   + 134\gamma_{3}\gamma_{9} - 123\gamma_{4}c_{8} 
                                 -674\gamma_{3}^{2}\gamma_{6} - 83\gamma_{3}\gamma_{4}\gamma_{5} )  \\
                       &  +t^{4}(139\gamma_{5}\gamma_{6} + 31\gamma_{4}c_{7} + 26\gamma_{3}c_{8} 
                          +117\gamma_{3}\gamma_{4}^{2} +130\gamma_{3}^{2}\gamma_{5}) \\
                       &  +t^{5}(513\gamma_{4}\gamma_{6}   -194\gamma_{3}c_{7} +604\gamma_{3}^{2}\gamma_{4} ) 
                          +t^{6}(\gamma_{9}+  1094\gamma_{3}\gamma_{6} + 133\gamma_{4}\gamma_{5})  
 \end{align*} 
 \begin{align*} 
                       &  +t^{7}(3c_{8} + 198\gamma_{3}\gamma_{5})   
                          +t^{8}(22c_{7}  -685\gamma_{3}\gamma_{4}) + t^{9}(\gamma_{6} + 18\gamma_{3}^{2})   \\
                       &  + 4 t^{10} \gamma_{5} + 241 t^{11} \gamma_{4}   + 382t^{12} \gamma_{3}. 
\end{align*}  
In particular, the set $\{ Z_{1}, Z_{2}, \ldots, Z_{8}, Z_{542}, Z_{6542}, Z_{76542}, Z_{136542}, Z_{154376542}, 
Z_{1654376542}, Z_{134276543876542} \}$ is a minimal system of ring generators of $H^*(E_{8}/T;\Z)$ that consists of 
Schubert classes.
\end{prop}

\section{Chow rings of $\E_{l} \; (l = 6, 7, 8)$}   \label{sec:Chow} 
In this section,  we determine the Chow rings of the exceptional  groups $\E_{l} \; (l=6,7,8)$.  
As mentioned in the introduction, 
we have only to compute the quotient ring of $A(G/B)$ by the ideal 
generated by $A^1(G/B)$.  
Denote by $p^*: A(G/B) \longrightarrow A(G)$ the  projection onto the quotient induced by the natural projection 
$p:G \to G/B$.  It is known that the cycle map $A(G/B) \to H^{*}(G/B;\Z)$ which assigns to the Schubert 
variety $X_{w_0w}$ the Schubert class $Z_{w}$ is an  isomorphism of rings.
Therefore, the Chow ring $A(G)$ is isomorphic to the quotient ring of 
$H^*(G/B;\Z) = H^*(K/T;\Z)$ by the ideal generated by $H^2(K/T;\Z)$.  
 
By Theorem \ref{thm:E_6/T},  we easily  obtain that  
\begin{equation*} 
       H^*(E_{6}/T;\Z)/(t_{1}, \ldots, t_{6}, t)  
             =  \Z[\gamma_{3}, \gamma_{4}]/(2\gamma_{3}, 3\gamma_{4}, \gamma_{3}^2, \gamma_{4}^3).  
 \end{equation*}  
Furthermore, by Proposition \ref{prop:E_6/T}, we have  
     \[  Z_{542} \equiv \gamma_{3}, \;  Z_{6542} \equiv \gamma_{4},  \mod (t_{1}, \ldots, t_{6}, t).    \] 

In a similar manner, by Theorem \ref{thm:E_7/T},  we obtain 
\begin{equation*} 
 H^*(E_{7}/T;\Z)/(t_{1}, \ldots, t_{7}, t)   
                = \Z[\gamma_{3}, \gamma_{4}, \gamma_{5}, \gamma_{9}]
                   /(2\gamma_{3}, 3\gamma_{4}, 2\gamma_{5}, \gamma_{3}^2, 2\gamma_{9}, \gamma_{5}^2, 
                     \gamma_{4}^3, \gamma_{9}^2).    
\end{equation*} 
Furthermore, by Proposition \ref{prop:E_7/T}, we have 
  \[  Z_{542} \equiv \gamma_{3}, \;  Z_{6542} \equiv \gamma_{4}, \;  Z_{76542} \equiv \gamma_{5},  \; 
      Z_{654376542}  \equiv \gamma_{9}  \mod (t_{1}, \ldots, t_{7}, t).   \]

Finally, by Theorem \ref{thm:E_8/T}, we obtain
\begin{align*} 
        & H^*(E_{8}/T;\Z)/(t_{1}, \ldots, t_{8}, t)   \\ 
                 = \ & \Z[\gamma_{3}, \gamma_{4}, \gamma_{5}, \gamma_{6}, \gamma_{9}, \gamma_{10}, \gamma_{15}] \\
                     &  /(2\gamma_{3}, 3\gamma_{4}, 2\gamma_{5}, 5\gamma_6,  2\gamma_{9}, 
                      \gamma_{5}^{2} - 3\gamma_{10}, \gamma_{4}^3,
                     2\gamma_{15},  \gamma_{9}^2,   3\gamma_{10}^{2}, \gamma_{3}^{8}, 
                      \gamma_{15}^{2} + \gamma_{10}^{3} + 2\gamma_{6}^{5}).    
\end{align*} 
Furthermore, by Proposition \ref{prop:E_8/T}, we have 
       \begin{align*}  
                Z_{542} & \equiv \gamma_{3}, \; Z_{6542} \equiv \gamma_{4},  \;  Z_{76542} \equiv \gamma_{5}, 
                Z_{136542} \equiv \gamma_{6},  \;  Z_{154376542} \equiv \gamma_{9}, \\
                Z_{1654376542} & \equiv -\gamma_{10} + \gamma_{5}^{2}, \;
                Z_{134276543876542}  \equiv \gamma_{15} +  \gamma_{5}\gamma_{10} + \gamma_{3}^{2} \gamma_{9} + \gamma_{3}^{5} 
                \mod (t_{1}, \ldots, t_{8}, t),
        \end{align*} 
and the relations are calculated as follows:
\begin{align*}
\gamma_5^2-3\gamma_{10}     & \equiv  Z_{76542}^2 - 3(-Z_{1654376542}+Z_{76542}^2) \equiv 3 Z_{1654376542}, \\
         3\gamma_{10}^2     & \equiv    Z_{76542}^2(-Z_{1654376542}+Z_{76542}^2) \equiv Z_{76542}^4, \\
           Z_{136542}^5     & \equiv  \gamma_6^5 \equiv -12(\gamma_{15}^2+\gamma_{10}^3+2\gamma_{6}^5), \; 
           Z_{1654376542}^3   \equiv  -\gamma_{10}^3 \equiv -10(\gamma_{15}^2+\gamma_{10}^3+2\gamma_{6}^5), \\
            \gamma_{15}     & \equiv  Z_{134276543876542} +Z_{76542} ^3+Z_{542}^2 Z_{154376542} +Z_{542}^5, \\
     Z_{134276543876542}^2  & \equiv  \gamma_{15}^2 \equiv 15(\gamma_{15}^2+\gamma_{10}^3+2\gamma_{6}^5) \;
                           \hfill \mod (t_{1}, \ldots, t_{8}, t).
\end{align*}

Therefore we obtain our main result of this paper: 
\begin{thm}  \label{thm:ChowE_6E_7E_8} 
{\rm (1)} The Chow ring of $\mathrm{E_6}$ is given by 
  \[  A(\E_{6}) = \Z[X_{3}, X_{4}]/(2X_{3}, 3X_{4}, X_{3}^2, X_{4}^3),  \]
  where $X_{3}$ and $X_{4}$ are the pull-back images of the elements of $A(\E_{6}/B)$  defined by 
  the Schubert varieties $X_{w_{0}s_{5}s_{4}s_{2}}$ and $X_{w_{0}s_{6}s_{5}s_{4}s_{2}}$ respectively. 

{\rm (2)} The Chow ring of $\mathrm{E_7}$ is given by 
     \[   A(\E_{7})  =   \Z[X_{3}, X_{4}, X_{5}, X_{9}]
                              /(2X_{3}, 3X_{4}, 2X_{5}, X_{3}^2, 2X_{9}, X_{5}^2, X_{4}^3, X_{9}^2),    \] 
  where $X_{3}, X_{4}, X_{5}$,  and $X_{9}$ are the pull-back images of the elements of  $A(\E_{7}/B)$ 
  defined by the Schubert varieties  $X_{w_{0}s_{5}s_{4}s_{2}}$, $X_{w_{0}s_{6}s_{5}s_{4}s_{2}}$, 
  $X_{w_{0}s_{7}s_{6}s_{5}s_{4}s_{2}}$, and $X_{w_{0}s_{6}s_{5}s_{4}s_{3}s_{7}s_{6}s_{5}s_{4}s_{2}}$
  respectively.
  
{\rm (3)} The Chow ring of $\E_{8}$ is given by
  \begin{align*} 
         A(\E_{8})  & =  \Z[X_{3}, X_{4}, X_{5}, X_{6}, X_{9}, X_{10}, X_{15}]  \\
                    & \hspace{0.3cm}/(2X_{3}, 3X_{4}, 2X_{5}, 5X_6,  2X_{9}, 3X_{10}, X_{4}^3, 2X_{15},  X_{9}^2, X_{5}^4, 
                       X_{3}^8, X_6^5, X_{10}^3, X_{15}^2),   
  \end{align*}   
  where $X_{3}, X_{4}, X_{5}, X_6, X_{9}, X_{10}$, and $X_{15}$ are the pull-back images under $p^*$ of the elements of  $A(\E_{8}/B)$ 
  defined by the Schubert varieties  $X_{w_0s_{5}s_{4}s_{2}}, X_{w_0s_{6}s_{5}s_{4}s_{2}}$, 
  $X_{w_0s_{7}s_{6}s_{5}s_{4}s_{2}}, X_{w_0s_{1}s_{3}s_{6}s_{5}s_{4}s_{2}}$, \\
  $X_{w_0s_{1}s_{5}s_{4}s_{3}s_{7}s_{6}s_{5}s_{4}s_{2}}$,
  $X_{w_0s_{1}s_{6}s_{5}s_{4}s_{3}s_{7}s_{6}s_{5}s_{4}s_{2}}$, and
  $X_{w_0s_{1}s_{3}s_{4}s_{2}s_{7}s_{6}s_{5}s_{4}s_{3}s_{8}s_{7}s_{6}s_{5}s_{4}s_{2}}$  respectively.
\end{thm}

\begin{rem}\rm 
As a Corollary of Theorem \ref{thm:ChowE_6E_7E_8},
we can recover the {\em $p$-exceptional degrees} of $G=\E_{l} \ (l = 6, 7, 8)$ for each torsion prime $p$.
In \cite{Kac85},  Kac showed that the kernel of the characteristic homomorphism tensored 
by a prime field $\F_p$   is an ideal generated by a regular sequence.
Hence, there  exist  homogeneous polynomials 
$P_{1},\ldots,P_{l}  \in S(\hat{H})\otimes \F_p \cong H^*(BT;\F_p)$   such that
\[
     \mathcal{I}_p := \ker{c\otimes \F_p} =  (P_{1},\ldots,P_{l}).
 \]
Then he introduced the notion of the $p$-exceptional degrees which is a certain sub-sequence of the degrees of 
the {\em basic invariants} of $G$.  The basic invariants of $G$ are defined to be  homogeneous generators of 
$\ker{c \otimes \Q} = \left( (S^+(\hat{H})\otimes \Q)^W \right)$.
We only show how to deal with the case of $(G,p)=(\E_{8},2)$ (cf. \cite[\S 7]{Kac85}). 
The other cases can be obtained in a similar manner.
By Theorem \ref{thm:ChowE_6E_7E_8},  the mod $2$ Chow ring of $\E_{8}$ is 
\begin{align*} 
     A(\E_{8}; \F_{2}) & = A(\E_{8}) \otimes_{\Z}  \F_{2}   \\
                      & = \F_{2}[X_{3}, X_{5}, X_{9}, X_{15}]/(X_{3}^{8}, X_{5}^{4}, X_{9}^{2}, X_{15}^{2}).    
\end{align*}  
Therefore,  by  \cite[Theorem 3]{Kac85}, 
the $2$-exceptional degrees are  $2 ^{1} \cdot 9 = 18$, $2^{2} \cdot 5 = 20$, $2^{3}  \cdot 3  = 24$, 
$2^{1} \cdot 15  = 30$.  Moreover, by \cite[Theorem 4]{Kac85}, 
 the degrees of generators of $\mathcal{I}_2$ can be read off  as follows:
in the sequence $(2, 8, 12, 14, 18, 20, 24, 30)$ of the degrees of the basic invariants  (see the table below),  
we replace the number $18$ by $9$, $20$ by $5$,  $24$ by $3$, and $30$ by $15$.

The $p$-exceptional degrees of $\E_{l} \ (l=6,7, 8)$ are given  by  the following table
(cf. \cite[Table 2]{Kac85}).

\begin{tabular}{l|lll}
$(G,p)$ & degrees of basic invariants & degrees of generators of $\mathcal{I}_p$ & $p$-exceptional degrees  \\
\hline
$(\mathrm{E_6},2)$ & $(2,5,6,8,9,12)$ & $(2,3,5,8,9,12)$ & $(6)$ \\
$(\mathrm{E_6},3)$ & $\uparrow$ & $(2,4,5,6,8,9)$ & $(12)$ \\ 
$(\mathrm{E_7},2)$ & $(2,6,8,10,12,14,18)$ & $(2,3,5,8,9,12,14)$ & $(6,10,18)$ \\
$(\mathrm{E_7},3)$ & $\uparrow$ & $(2,4,6,8,10,14,18)$ & $(12)$ \\
$(\E_8,2)$ & $(2,8,12,14,18,20,24,30)$ & $(2,3,5,8,9,12,14,15)$ & $(18,20,24,30)$ \\
$(\E_8,3)$ & $\uparrow$ & $(2,4,8,10,14,18,20,24)$ & $(12,30)$\\
$(\E_8,5)$ & $\uparrow$ & $(2,6,8,12,14,18,20,24)$ & $(30)$
\end{tabular}
\end{rem}

\section{Appendix} 
In this appendix, 
we first give the Borel presentations  of the integral cohomology rings of $E_{l}/T \ (l = 6, 7, 8)$
in terms of the rings of invariants of the Weyl groups $W(E_{l})$.  
Then we  relate these results  to the Schubert presentations 
obtained by  Duan-Zhao (\cite{DZ08}).

\subsection{Borel presentation of $H^*(E_{l}/T;\Z)$}  \label{Borel_presentation}  
As mentioned in \S \ref{subsec:E_l/T}, we derive $H^*(E_{l}/T;\Z)$ directly 
from the ring of invariants of the Weyl group $W(E_{l})$.  
The notation used here is the same as in \S \ref{subsec:E_l/T}.

First we recall   the ring of invariants of the Weyl group $W(E_{l})$  in \cite[\S 5]{Toda-Wat74},
 \cite[\S 2]{Wat75},  and \cite[\S 2]{Nak09} (see also \cite[2.1, 2.2, 2.3]{Meh88}).  
 Putting 
  \begin{equation}  \label{eqn:x_i}
     \begin{array}{llll}  
       &  x_{i} = 2t_{i}-t \;  (1 \leq i \leq 6)  \quad (l = 6),   \medskip   \\
       &  x_{i} = 2t_{i} - t \; (1 \leq i \leq 7) \quad \text{and} \quad x_{8} = t  \quad (l = 7),  \medskip   \\
       &  x_{i} = 2t_{i} - \dfrac{2}{3}t \; (1 \leq i \leq 8)  \quad \text{and} \quad x_{9}= -\dfrac{2}{3}t   \quad (l = 8),  
     \end{array} 
 \end{equation}    we see  easily that the sets 
 \begin{align*} 
       &  S_{6} =   \{ x_{i}+x_{j} \; (1 \leq i < j \leq 6), \;  t-x_{i},-t-x_{i} \;    (1 \leq i \leq 6)  \} \quad (l = 6),   \\
       &  S_{7} =   \{ x_{i}+x_{j},  -x_{i}-x_{j} \; (1 \leq i < j \leq 8) \} \quad (l = 7), \\
       &  S_{8} =   \{ \pm \;(x_{i} - x_{j}) \; (1 \leq i < j \leq 9), \; \pm \; (x_{i} + x_{j} + x_{k}) \; 
                    (1 \leq i < j < k \leq 9)  \}  \quad (l = 8)
 \end{align*} 
are  invariant under the action of $W(E_{l}) \ (l = 6, 7, 8)$ respectively. 
Then  we form   the $W(E_{l})$-invariant polynomials 
    \[   I_{n} = \sum_{y \in S_{l}} y^{n} \in H^{2n}(BT;\Q)^{W(E_{l})}.  \]
The invariant polynomials  $I_{n}$  are  computed by the formula: 
  \begin{equation} \label{eqn:inv.forms.E_6E_7}  
   \begin{array}{llll} 
      I_{n} & \hspace{-0.2cm} = \dfrac{1}{2}\displaystyle{\sum_{i=2}^{n-2}} \binom{n}{i}s_{i}s_{n-i}  + (6-2^{n-1})s_{n} 
                + 2(-1)^{n}\sum_{j=0}^{[\frac{n}{2}]}  \binom{n}{2j}s_{n-2j}t^{2j}  \quad (l = 6),  \medskip \\    
      I_{n} & \hspace{-0.2cm} = (16-2^{n})s_{n} + \displaystyle{\sum_{i = 2}^{n  - 2} }
              \binom{n}{i}s_{i}s_{n-i}  \quad \text{ for } \; n  \;    \;  \text{even}  \quad (l = 7),  \medskip  \\
       I_{n} & \hspace{-0.2cm} =  \displaystyle{\sum_{i = 0}^{n}} \binom{n}{i}(-1)^{n-i} s_{i}s_{n-i} 
              + 2 \cdot 3^{n-1} s_{n}    - \displaystyle{\sum_{i = 0}^{n}}   \binom{n}{i}2^{n-i}s_{i}s_{n-i}   \medskip \\
            &  + \dfrac{1}{3}  \displaystyle{\sum_{i = 0}^{n} \sum_{j = 0}^{n-i}}
                 \binom{n}{i}  \binom{n-i}{j} s_{i}s_{j}s_{n-i-j}    
               \quad \text{for} \; n \; \; \text{even} \quad (l = 8), 
    \end{array} 
 \end{equation}   
where $s_{n} = x_{1}^{n} + \cdots + x_{6}^{n} \;  (l = 6)$, $s_{n} = x_{1}^{n} +  \cdots + x_{8}^{n} \;  (l = 7)$, 
      $s_{n} = x_{1}^{n} + \cdots + x_{9}^{n} \;  (l = 8)$.  
The power sum symmetric polynomials $s_{n}$ are written as polynomials in $d_{i}$'s, where  
$d_{i} = e_{i}(x_{1},\ldots,x_{6}) \;  (l = 6)$,
$d_{i} = e_{i}(x_{1}, \ldots, x_{8}) \;  (l = 7)$,  $d_{i} = e_{i}(x_{1}, \ldots, x_{9}) \;  (l = 8)$ by  use
of the Newton formula:
 \begin{equation}  \label{eqn:newton.formula}
     s_{n} = \sum_{i=1}^{n-1}(-1)^{i-1}s_{n-i}d_{i} + (-1)^{n-1}nd_{n}.   
  \end{equation}
Moreover, using (\ref{eqn:x_i}), we have 
   \begin{equation}  \label{eqn:kakikae}
     \begin{array}{lll} 
      d_{n}  & \hspace{-0.2cm} = \displaystyle{\sum_{i=0}^{n} }\binom{6-i}{n-i}(-t)^{n-i}2^{i}c_{i} \quad (l = 6),   \medskip \\
      d_{n}  & \hspace{-0.2cm} = \displaystyle{\sum_{i=0}^{n}  }  \left \{ \binom{7-i}{n-i} 
                                             - \sum_{i = 0}^{n}  \binom{7 - i}{n - 1 - i} \right \}   
                                             (-t)^{n - i} 2^{i}c_{i}  \quad (l = 7), \medskip \\
      d_{n}  & \hspace{-0.2cm} =  \displaystyle{\sum_{i=0}^{n} }  \left \{ \binom{8-i}{n-i}  +  \binom{8-i}{n--i} \right \}
                                    \left (-\frac{2}{3}t \right )^{n-i}2^{i}c_{i}  \quad (l = 8). \medskip      
     \end{array}  
 \end{equation}        
By using (\ref{eqn:inv.forms.E_6E_7}), (\ref{eqn:newton.formula}),  and (\ref{eqn:kakikae}), 
 $I_{n}$ can be written as  polynomials in $t$ and $c_{2}, \ldots, c_{l}$.  
Then the next lemma is proved in \cite{Toda-Wat74}, \cite{Wat75}, and \cite{Nak09} (see also \cite{Meh88}). 
\begin{lem}[\cite{Toda-Wat74}, Lemma 5.2, \cite{Wat75}, Lemma 2.1, \cite{Nak09}, Lemma 2.6]
  The rings  of invariants of  the Weyl group $W(E_{l}) \ (l = 6, 7, 8)$  with rational coefficients are  
   respectively given as follows$:$ 
   \begin{align*} 
          H^{*}(BT;\Q)^{W(E_{6})} & = \Q[I_{2}, I_{5}, I_{6}, I_{8},I_{9},I_{12} ],  \\
          H^{*}(BT;\Q)^{W(E_{7})} & = \Q[I_{2}, I_{6}, I_{8}, I_{10},I_{12},I_{14}, I_{18}], \\
          H^{*}(BT;\Q)^{W(E_{8})} & = \Q[I_{2}, I_{8}, I_{12}, I_{14}, I_{18}, I_{20}, I_{24}, I_{30}]. \\        
   \end{align*} 
\end{lem}
Hence, by Theorem \ref{thm:Borel}, the rational cohomology rings of $E_{l}/T \ (l = 6, 7, 8)$ are  
   \begin{align*} 
         H^*(E_{6}/T;\Q)  &=  H^*(BT;\Q)/(H^{+}(BT;\Q)^{W(E_{6})})  \\ 
                          &=  \Q[t_{1}, t_{2}, \ldots, t_{6}]/(I_{2}, I_{5}, I_{6}, I_{8}, I_{9}, I_{12}),   \\
         H^*(E_{7}/T;\Q)  &=  H^*(BT;\Q)/(H^{+}(BT;\Q)^{W(E_{7})})  \\ 
                          &=  \Q[t_{1}, t_{2}, \ldots, t_{7}]/(I_{2}, I_{6}, I_{8}, I_{10}, I_{12}, I_{14}, I_{18}), \\              
         H^*(E_{8}/T;\Q)  &= H^*(BT;\Q)/(H^{+}(BT;\Q)^{W(E_{8})})  \\ 
                          &=  \Q[t_{1}, t_{2}, \ldots, t_{8}]/(I_{2}, I_{8}, I_{12}, I_{14}, I_{18}, I_{20}, I_{24}, I_{30}),
   \end{align*} 
Then, by \cite[Theorem 2.1]{Toda75},  the integral relations $\rho_{j}$ 
are determined by the maximal integers $n_{j}$ in 
     \[  n_{j}  \cdot \rho_{j} \equiv I_{j}  \mod (\rho_{i}; i < j).   \]   
These integers are given by 
\begin{center}  
   \begin{tabular}{lllllllll} 
   \noalign{\hrule height0.8pt}
      $E_{6}$                                  & &  $E_{7}$                                              & &  $E_{8}$    \\
   \hline 
      $n_{2}  = -2^{4} \cdot 3$                & &  $n_{2}  = -2^{5} \cdot 3$                            & &  $n_{2}  = -2^{5} \cdot 3 \cdot 5$   \\
      $n_{5}  = -2^{7} \cdot 3 \cdot 5$        & &  $n_{6}  = 2^{10} \cdot 3^{2}$                        & &  $n_{8}  = 2^{15} \cdot 3^{2} \cdot 5$  \\
      $n_{6}  = 2^{9} \cdot 3^{2}$             & &  $n_{8}  = 2^{13} \cdot 3 \cdot 5$                    & &  $n_{12} = 2^{18} \cdot 3^{4} \cdot 5 \cdot 7$ \\
      $n_{8}  = 2^{12} \cdot 3 \cdot 5$        & &  $n_{10} = 2^{14} \cdot 3^{2} \cdot 5 \cdot 7$        & &  $n_{14} = 2^{20} \cdot 3^{2} \cdot 5^{2} \cdot 7 \cdot 11$  \\
      $n_{9}  = 2^{11} \cdot 3^{3} \cdot 7$    & &  $n_{12} = -2^{16} \cdot 3^{4} \cdot 5$               & &  $n_{18} = 2^{26} \cdot 3^{4} \cdot 5^{2} \cdot 7 \cdot 13$  \\
      $n_{12} = -2^{15} \cdot 3^{4} \cdot 5$   & &  $n_{14} = 2^{17} \cdot 3 \cdot 7 \cdot 11 \cdot 29$  & &  $n_{20} = 2^{27} \cdot 3^{2} \cdot 5^{2} \cdot 11 \cdot 17 \cdot 41$  \\
                                               & &  $n_{18} = 2^{22} \cdot 3^{3} \cdot 1229$             & &  $n_{24} = 2^{32} \cdot 3^{3} \cdot 5 \cdot 7 \cdot 11 \cdot 19 \cdot 199$  \\
                                               & &                                                       & &  $n_{30} = 2^{37} \cdot 3^{4} \cdot 5^{5} \cdot 7 \cdot 11 \cdot 13 \cdot 61$  \\
 \noalign{\hrule height0.8pt}
   \end{tabular} 
\end{center}

\subsection{Schubert presentation of $H^*(E_{l}/T;\Z)$} 
For the determination of $A(\E_{l})$, it is enough to give a  minimal system of generators of 
$H^*(E_{l}/T;\Z)$   that consists of Schubert classes as we  did  in \S \ref{subsec:BGG}. 
However, it would be  desirable to  write   down the relations explicitly among   these Schubert classes, and 
give  the presentation of $H^*(E_{l}/T;\Z)$ in terms of Schubert classes. 
Recently, one such  presentation has been   obtained  by Duan-Zhao \cite{DZ08}  in a different manner.  
In this subsection, we compare the generators and the relations in Theorems \ref{thm:E_6/T}, \ref{thm:E_7/T}, 
and \ref{thm:E_8/T}  with  those given  in \cite{DZ08}.

Following Duan-Zhao \cite{DZ08}, we set 
 \begin{equation}  \label{eqn:gen.(DZ)}  
   \begin{array}{llll}  
        y_{3} & = Z_{542},   \;   y_{4} = Z_{6542} \; (l = 6, 7, 8),  \medskip \\
        y_{5} & = Z_{76542},    \;   y_{9} = Z_{154376542} \; (l = 7, 8).  \medskip  \\
        y_{6} & = Z_{136542},    \;   y_{10} = Z_{1654376542}, \;  y_{15} = Z_{542316543876542}  \; (l = 8).  \medskip         
  \end{array}
\end{equation}   
In \cite{DZ08},  the fundamental weights $\{ \omega_{i} \}_{1 \leq i \leq l} \ (l = 6, 7, 8)$ are regarded 
as elements of $H^{2}(E_{l}/T;\Z) \ (l = 6, 7, 8)$ via the isomorphism 
$c: H^{2}(BT;\Z) \longrightarrow H^{2}(E_{l}/T;\Z)$.  Therefore 
$Z_{i} = \omega_{i} \; (1 \leq i  \leq l) \ (l = 6, 7, 8)$ and $t = \omega_{2}$ in their notation.

Then Duan-Zhao obtained the following description of the integral cohomology ring of $E_{6}/T$: 
\begin{thm}[Duan-Zhao \cite{DZ08}, Theorem 3]  \label{thm:E_6/T(DZ)} 
    \[  H^*(E_{6}/T;\Z) = \Z[\omega_{1}, \ldots, \omega_{6}, y_{3}, y_{4}]
                                                        /(r_{2}, r_{3}, r_{4}, r_{5}, r_{6}, r_{8}, r_{9}, r_{12}),  \] 
where 
\begin{align*} 
  r_{2} &= 4t^2 - c_{2}, \;    r_{3} = 2y_{3} + 2t^3 - c_{3}, \;   r_{4} = 3y_{4} + t^4 - c_{4}, \\ 
  r_{5} &= 2t^{2}y_{3} - tc_{4} + c_{5}, \;   r_{6} = y_{3}^2 - tc_{5} + 2c_{6}, \\
  r_{8} &= 3y_{4}^2 - 2c_{5}y_{3} - t^{2}c_{6} + t^3 c_{5}, \;   r_{9} = 2y_{3}c_{6} - t^3 c_{6},  \;  r_{12} = y_{4}^3 - c_{6}^2.  
\end{align*} 
\end{thm}

Then, by  Theorems \ref{thm:E_6/T}, \ref{thm:E_6/T(DZ)},  Proposition \ref{prop:E_6/T},  and (\ref{eqn:gen.(DZ)}),  
the correspondence is given as follows:

\begin{prop}  
\begin{align*} 
 y_{3} & = \gamma_{3} - t^3, \;   y_{4} = \gamma_{4} - t^4;   \\
 r_{2} & = -\rho_{2}, \;   r_{3}  = -\rho_{3}, \;   r_{4}  = -\rho_{4}, \;  r_{5}  = -\rho_{5} + t \rho_{4}, \; 
 r_{6}  = \rho_{6} - t\rho_{5}, \\
 r_{8} & = \rho_{8} - 2y_{3}\rho_{5} + 4t^{2}\rho_{6} + t^{3} \rho_{5}, \;  r_{9}  = -\rho_{9}, \; 
 r_{12}  = \rho_{12} + y_{4}\rho_{8} + y_{3} \rho_{9} - 2c_{6} \rho_{6}. 
\end{align*} 

\end{prop}

Likeweise, their  description of the integral cohomology ring of $E_{7}/T$ is given as follows:  
\begin{thm} [Duan-Zhao \cite{DZ08}, Theorem 4] \label{thm:E_7/T(DZ)}
    \[     H^*(E_{7}/T;\Z)  =  \Z[\omega_{1}, \ldots, \omega_{7}, y_{3}, y_{4}, y_{5}, y_{9}] 
                        /(r_{2}, r_{3}, r_{4}, r_{5}, r_{6}, r_{8}, r_{9}, r_{10}, r_{12}, r_{14}, r_{18}), \] 
where 
\begin{align*} 
  r_{2} &= 4t^2 - c_{2}, \;     r_{3} = 2y_{3} + 2t^3 - c_{3}, \;  r_{4} = 3y_{4} + t^4 - c_{4}, \\
  r_{5} &= 2y_{5} - 2t^2 y_{3} + tc_{4} - c_{5}, \;   r_{6} = y_{3}^2 - tc_{5} + 2c_{6}, \\
  r_{8} &= 3y_{4}^2 + 2y_{3}y_{5} - 2y_{3}c_{5} + 2tc_{7} - t^2 c_{6} + t^3 c_{5}, \\
  r_{9} &= 2y_{9} + 2y_{4}y_{5} - 2y_{3}c_{6} - t^2 c_{7} + t^3 c_{6}, \;   
  r_{10} = y_{5}^2 - 2y_{3}c_{7} + t^3 c_{7}, 
\end{align*} 
\begin{align*} 
  r_{12} &= y_{4}^3 - 4y_{5}c_{7} - c_{6}^2 - 2y_{3}y_{9} - 2y_{3}y_{4}y_{5} + 2ty_{5}c_{6} + 3t y_{4}c_{7}
            + c_{5}c_{7}, \\
  r_{14} &= c_{7}^2 - 2y_{5}y_{9} + 2y_{3}y_{4}c_{7} - t^3 y_{4}c_{7}, \\
  r_{18} &= y_{9}^2 + 2y_{5}c_{6}c_{7} - y_{4}c_{7}^2 - 2y_{4}y_{5}y_{9} + 2y_{3}y_{5}^3 - 5ty_{5}^2 c_{7}. 
\end{align*}  
\end{thm} 

By Theorems \ref{thm:E_7/T}, \ref{thm:E_7/T(DZ)},  Proposition \ref{prop:E_7/T}, and (\ref{eqn:gen.(DZ)}),  
the correspondence is given as follows: 
\begin{prop} 
\begin{align*}
 y_{3} & = \gamma_{3} -  t^3, \; y_{4} = \gamma_{4} -  t^4, \; y_{5} = \gamma_{5}, \; 
 y_{9} = \gamma_{9} - \gamma_{4}\gamma_{5} +t^{4} \gamma_{5}  \\
 r_{2} &= -\rho_{2}, \;  r_{3} = -\rho_{3}, \;  r_{4} = -\rho_{4}, \;  r_{5} = -\rho_{5} + t \rho_{4}, \;   
 r_{6} =  \rho_{6} - t\rho_{5}, \\  
 r_{8} &= \rho_{8} - 2y_{3}\rho_{5} + 4t^2 \rho_{6} + t^3 \rho_{5}, \;  r_{9} = -\rho_{9}, \;  r_{10} = \rho_{10}, \\
 r_{12} &= \rho_{12} + y_{4}\rho_{8} + y_{3}\rho_{9} - 2c_{6}\rho_{6} + c_{7}\rho_{5}, \; 
 r_{14} = \rho_{14} + y_{5}\rho_{9} + 2y_{4}\rho_{10}, \\
 r_{18} &= \rho_{18} + y_{4}\rho_{14} + (3y_{4}^2 + 2y_{3}y_{5} - 5tc_{7})\rho_{10} + (y_{4}y_{5} - y_{9})\rho_{9} 
           + (-2c_{7}y_{3} + t^3 c_{7})\rho_{8}  \\
         & + (-12tc_{7}y_{4}- 24t^5 c_{7})\rho_{6}. 
\end{align*} 
\end{prop}

For the case of $E_{8}/T$,  they gave the following presentation:   
\begin{thm} [Duan-Zhao \cite{DZ08}, Theorem 5]  \label{thm:E_8/T(DZ)}  
\begin{align*} 
        H^*(E_{8}/T;\Z)  = \ & \Z[\omega_{1}, \ldots, \omega_{8}, y_{3}, y_{4}, y_{5}, y_{6}, y_{9}, y_{10}, y_{15}]  \\
                             & /(r_{2}, r_{3}, r_{4}, r_{5}, r_{6}, r_{8}, r_{9}, r_{10}, r_{12}, r_{14}, r_{15}, r_{18}, 
                             r_{20}, r_{24}, r_{30}), 
\end{align*} 

where 
\begin{align*} 
r_{2} &= 4\omega _{2}^{2}-c_{2}, \;
r_{3} = 2y_{3}+2\omega _{2}^{3}-c_{3}, \;
r_{4}=3y_{4}+\omega _{2}^{4}-c_{4},  \\
r_{5} &=2y_{5}-2\omega _{2}^{2}y_{3}+\omega _{2}c_{4}-c_{5}, \;
r_{6}=5y_{6}+2y_{3}^{2}+10\omega _{2}y_{5}-2\omega _{2}c_{5}-c_{6}, \\
r_{8} &=3{c_{8}}-3{y_{4}^{2}}-2{y_{3}y_{5}}+2{y_{3}}c_{5}-{2\omega _{2}c_{7}}+
{\omega _{2}^{2}c_{6}-\omega _{2}^{3}c_{5}}, \\
r_{9}&=2{y_{9}}+2{y_{4}y_{5}}-2{y_{3}y_{6}}-4{{\omega _{2}}y_{3}y_{5}+
{\omega _{2}c_{8}}-\omega _{2}^{2}c_{7}}+{\omega _{2}^{3}c_{6}}, \\
r_{10}&=3{{y_{10}-2{y_{5}^{2}}-2{y_{3}c_{7}}-3{y_{4}y_{6}+3{y_{4}}c_{6}}
 -6\omega _{2}{y_{4}y_{5}}}}-{\omega _{2}^{2}c_{8}}+{\omega _{2}^{3}c_{7}}, \\
r_{12}\mid _{\omega_{2}=0} &= y_{4}^{3}-2c_{4}c_{8}-c_{5}c_{7}+3c_{6}y_{6}+c_{3}^{2}y_{6}+c_{3}y_{3}^{3}, \\
r_{14}\mid _{\omega_{2}=0}&=c_{7}^{2}+c_{4}y_{10}-c_{3}^{2}c_{8}-c_{4}y_{4}y_{6}, \\
r_{15}\mid _{\omega_{2}=0} &=2y_{15}+c_{5}y_{10}+5c_{7}c_{8}-c_{8}y_{3}y_{4}-2c_{3}c_{5}c_{7}+c_{3}c_{6}y_{6}+2c_{5}y_{4}y_{6} -c_{3}y_{3}^{2}y_{6}+c_{4}y_{3}y_{4}^{2}+c_{3}y_{4}^{3}, \\
r_{18}\mid _{\omega_{2}=0} &=y_{9}^{2}-6y_{10}y_{8}-4y_{9}y_{6}y_{3}+4y_{9}y_{5}y_{4}+5y_{8}y_{4}y_{3}^{2}+y_{7}^{2}y_{4}
-3y_{7}y_{4}^{2}y_{3}+5y_{6}^{3}+3y_{6}^{2}y_{3}^{2}+10y_{6}y_{5}y_{4}y_{3}  \\
&+y_{6}y_{4}^{3}+6y_{5}y_{4}y_{3}^{3}, \\
r_{20}\mid _{\omega_{2}=0}&=(y_{10}+4y_{4}y_{6}-y_{5}^{2}+2y_{3}^{2}y_{4})^{2}-y_{8}(6y_{3}y_{9}+3y_{4}y_{8}-y_{5}y_{7} +14y_{3}^{2}y_{6}+8y_{3}^{4}), \\
r_{24}\mid _{\omega_{2}=0}&=5(y_{3}^{2}+2y_{6})^{4}-y_{8}(5y_{7}y_{9}+5y_{8}^{2}-4y_{3}y_{5}y_{8}+2y_{3}^{3}y_{7}+20y_{5}^{2}y_{6} 
+10y_{3}^{2}y_{4}y_{6}+18y_{3}y_{4}^{2}y_{5}+4y_{3}^{4}y_{4}), \\
r_{30}\mid _{\omega_{2}=0}&=(y_{3}^{2}+2y_{6})^{5}+(y_{10}+4y_{6}y_{4}-y_{5}^{2}+2y_{4}y_{3}^{2})^{3}+(y_{15}+y_{10}y_{5} 
+y_{9}y_{3}^{2}+2y_{8}y_{7}-4y_{7}y_{5}y_{3}+5y_{6}^{2}y_{3} \\
&+2y_{6}y_{5}y_{4}+2y_{5}y_{4}y_{3}^{2}+y_{4}^{3}y_{3})^{2}
+y_{8}(y_{15}y_{7}+8y_{15}y_{4}y_{3}-9y_{10}y_{8}y_{4}-10y_{10}y_{6}y_{3}^{2}-4y_{10}y_{4}^{3}-2y_{10}y_{3}^{4} \\
&+y_{9}y_{7}y_{3}^{2}+6y_{9}y_{5}y_{4}^{2}+8y_{9}y_{4}y_{3}^{3}-2y_{8}^{2}y_{6}-y_{8}^{2}y_{3}^{2}+44y_{8}y_{7}^{2}+7y_{8}y_{6}y_{5}y_{3} 
-49y_{8}y_{5}^{2}y_{4}+7y_{8}y_{5}y_{3}^{3} \\
&+25y_{7}^{2}y_{5}y_{3}-5y_{7}y_{6}^{2}y_{3}+10y_{7}y_{6}y_{5}y_{4} 
-12y_{7}y_{5}y_{4}y_{3}^{2}-30y_{7}y_{4}^{3}y_{3}-10y_{6}^{3}y_{4}+5y_{6}^{2}y_{5}^{2}+12y_{5}^{3}y_{4}y_{3}  \\&+3y_{5}y_{4}^{2}y_{3}^{3} 
+y_{4}^{4}y_{3}^{2}+4y_{4}y_{3}^{6}).
\end{align*}  
\end{thm}

By Theorems \ref{thm:E_8/T}, \ref{thm:E_8/T(DZ)}, Proposition \ref{prop:E_8/T}, and (\ref{eqn:gen.(DZ)}), the correspondence  is given as follows:
\begin{prop}
 \begin{align*}
 y_{3} & = \gamma_{3} - t^{3}, \; y_{4} = \gamma_{4} - t^{4},  \;  y_{5} = \gamma_{5},  
 y_{6} = \gamma_{6} - t\gamma_{5} + t^{2} \gamma_{4}, \\
 y_{9} &= \gamma_{9} -  2\gamma_{3}^{3} - 4\gamma_{3}\gamma_{6}        
                     - \gamma_{4} \gamma_{5} + t(- 6\gamma_{4}^{2} + 5c_{8} - 4 \gamma_{3}\gamma_{5}) 
                     + t^{2} (- 4c_{7} + 14\gamma_{3} \gamma_{4}) 
                     + t^{3} (- 2\gamma_{3}^{2} + 14\gamma_{6}) \\
                    & - 5 t^{4} \gamma_{5}  - 10 t^{5} \gamma_{4} + 10t^{6}\gamma_{3}, \\
  y_{10} &= -\gamma_{10} +  \gamma_{5}^{2}  -2\gamma_{3}^{2}\gamma_{4} -4\gamma_{4}\gamma_{6} 
                     + 2t^{2}\gamma_{4}^{2} +  t^{4}(2\gamma_{3}^{2} + 4\gamma_{6}) - 4t^{6}\gamma_{4} + 2t^{10}, \\
  y_{15} & \equiv  \gamma_{15}   +   10\gamma_{6}\gamma_{9} -  2\gamma_{5}\gamma_{10} +   5\gamma_{3}\gamma_{6}^{2}
                      -  \gamma_{3}^{3} \gamma_{6} +  4\gamma_{3}^{2} \gamma_{9} + 2\gamma_{4}\gamma_{5}\gamma_{6}   
                      + 3\gamma_{3}^{2} \gamma_{4}\gamma_{5}  \mod (\omega_{2}).  
 \end{align*}
\end{prop}


\end{document}